\documentclass[12pt]{amsart}
\usepackage[usenames]{color} 
\usepackage{times}
\usepackage{amsfonts}
\usepackage{amsthm}
\usepackage{amsmath}
\usepackage{psfrag}

\usepackage{times}
\usepackage{latexsym,color}
\usepackage{amssymb}
\usepackage{graphicx}

\usepackage{epsfig}

\usepackage{gastex}

\usepackage{tikz}
\usepackage[all]{xy}
 \usepackage[hang]{caption}
 \usetikzlibrary{snakes}
\usetikzlibrary{arrows}

\advance \topmargin by -\headheight
\advance \topmargin by -\headsep
\evensidemargin \oddsidemargin
\def \ind{_{n \in {\mbox{\rm {\scriptsize I$\!$N}}}}}

\newcommand{\GZ}{{\mathbb Z}}

\newcommand{\GN}{{\mathbb N}}

\newcommand{\ab}{|}

\newtheorem{theorem}{Theorem}
\newtheorem{lemma}[theorem]{Lemma}
\newtheorem{corollary}[theorem]{Corollary}
\newtheorem{pp}[theorem]{Proposition}

\newtheorem{definition}{Definition}

\begin{document}

\title[Veech 1969]{Rigidity of generalized Veech 1969/Sataev 1975   extensions of rotations}

\author[S. Ferenczi]{S\'ebastien Ferenczi} 
\address{Aix Marseille Universit\'e, CNRS, Centrale Marseille, Institut de Math\' ematiques de Marseille, I2M - UMR 7373\\13453 Marseille, France.}
\email{ssferenczi@gmail.com}
\author[P. Hubert]{Pascal Hubert} 
\address{Aix Marseille Universit\'e, CNRS, Centrale Marseille, Institut de Math\' ematiques de Marseille, I2M - UMR 7373\\13453 Marseille, France.}
\email{hubert.pascal@gmail.com}

\date{June 30, 2020}

\begin{abstract}
We look at $d$-point extensions of a rotation of angle $\alpha$  with $r$ marked points, generalizing the examples of Veech 1969 and Sataev 1975, together with the square-tiled interval exchange transformations of \cite{fh2}. We study the property of rigidity, in function of the Ostrowski expansions of the marked points by $\alpha$: we  prove  that $T$ is rigid  when $\alpha$ has unbounded partial quotients, and that $T$ is  not rigid when the natural coding of the underlying rotation with  marked points is linearly recurrent. But there remains an interesting grey zone between these two cases, in which we have only partial results on the rigidity question; they allow us to build the first examples of non linearly recurrent and non rigid interval exchange transformations.

\end{abstract}
\maketitle

In a founding paper of 1969 \cite{ve69} W.A. Veech defines an extension of a rotation of angle $\alpha$ to two copies of the torus with a marked point $\beta$, the change of copy occurring on the interval $[0,\beta[$ (resp. $[\beta,1[$ on a variant): for $\alpha$ with unbounded partial quotients and some values of  $\beta$, they provide examples of minimal non uniquely ergodic interval exchange transformations. These systems were defined again independently, in a generalized way,   by E.A. Sataev in 1975, in a beautiful but not very well known paper \cite{sat}: by taking $r$ {\em marked points} and $r+1$ copies of the torus, he gets minimal interval exchange transformations with a prescribed number of ergodic invariant measures; also, improvements on Veech's results were introduced by  M. Stewart  \cite{ste} and K.D. Merrill \cite{mer}.  In the present paper, we study more general systems, by marking $r$ points and taking $d$ copies of the torus, for any $r\geq 1$, $d\geq 2$, and extending the rotation by the symmetric group $S_{d}$.

Though in general  our marked points are not in $\mathbb Z(\alpha)$, we allow one of them to be $1-\alpha$, so that our systems generalize also the square-tiled interval exchange transformations studied in \cite{fh2}. In that paper, we focussed on the measure-theoretic property of rigidity, meaning that for some sequence $q_n$ the $q_n$-th powers of the transformation  converge to the identity (Definition \ref{drig} below). 
 Experimentally, in the class of interval exchange transformations, the absence of rigidity is difficult to achieve (indeed, by Veech \cite{vs} it is true only for a set of measure zero of parameters) and all known examples satisfy also the word-combinatorial property of {\em linear recurrence} (Definition \ref{dlr} below) for their natural coding. Indeed, for the systems studied in \cite{fh2}, we proved rigidity is equivalent to absence of  linear recurrence, and thus to $\alpha$ having unbounded partial quotients. For the more general systems considered in the present paper, adaptations of the techniques of \cite{fh2} do allow us to prove
  that $T$ is rigid (for every invariant measure) and not linearly recurrent when $\alpha$ has unbounded partial quotients, and that $T$ is  uniquely ergodic, linearly recurrent, and not rigid when the natural coding of the underlying rotation with  marked points is linearly recurrent (under an extra  condition on the permutations, we prove also that $T$ is not of rank one); this linear recurrence requires $\alpha$ to have bounded partial quotients and the marked points to satisfy some conditions on their Ostrowski expansions by $\alpha$.

But, in sharp contrast with \cite{fh2}, there is no similar equivalence in the present class of systems, and this leaves an interesting grey zone, when $\alpha$ has bounded partial quotients but the Ostrowski expansions of the marked points do not satisfy the conditions required for linear recurrence; this means that   in the {\em Rokhlin towers} defined by the rotation (see Section \ref{roto} below) the marked points $\beta_i$ come too close either to one another or to the points $0$, $\alpha$, $1-\alpha$.  In these cases we prove some partial results, namely a sufficient condition (Theorem \ref{nrig}) for non-rigidity (for every ergodic invariant measure) and a sufficient condition (Theorem \ref{tcom}) for rigidity (for every  invariant measure): the latter puts all the grey zone on the rigid side for Veech 1969, while, with two (or more) marked points, the former  allows us to build the first known examples of non linearly recurrent and non rigid interval exchange transformations, answering Question 8 of \cite{fh2}:

\begin{theorem}\label{nrnlr} There exists a two-point extension of a rotation with two marked points which is a non linearly recurrent and non rigid interval exchange. \end{theorem}

These conditions are enough to give a full characterization of rigidity for the simplest generalizations of Veech 1969, when we take two copies of the torus and a small number of marked points (for higher numbers of marked points, the question is not untractable but  the results become extremely tedious to state). In general, the grey zone seems quite complicated, with many different cases using different  techniques, and we seem to be far from a complete characterization of rigidity in our class.

\section{Definitions} \label{sec:def}

\subsection{Word combinatorics}
We begin with basic definitions.
We look at finite {\em words} on a finite alphabet ${\mathcal A}=\{1,...k\}$. A word $w_1...w_s$ has
{\em length} $\ab w\ab=s$ (not to be confused with the length of a corresponding interval).   The {\em concatenation} of two words $w$ and $w'$ is denoted by $ww'$. 

\begin{definition}\label{dln} 
\noindent A word $w=w_1...w_s$ {\em occurs at place $i$} in a word $v=v_1...v_{s'}$ or an infinite sequence  $v=v_1v_2...$ if $w_1=v_i$, ...$w_t=v_{i+s-1}$. We say that $w$ is a {\em factor} of $v$. \\
A {\em   language} $L$ over $\mathcal A$ is a set of words such if $w$ is in $L$, all its factors are in $L$,  
\noindent A language $L$ is
{\em  minimal} if for each $w$ in $L$ there exists $n$ such that $w$
occurs in each word  of $L$ with $n$ letters.

\noindent The language $L(u)$ of an infinite sequence $u$ is the set of its finite factors.\\
A word $w$ is called {\em  right special}, resp. {\em
left
special} if there are at least two different letters $x$ such that $wx$, resp. $xw$, is in $L$. If $w$ is both right special and
left special, then $w$ is called {\em  bispecial}. \end{definition}

\subsection{Symbolic dynamics and codings}\label{sdc}

\begin{definition} The {\em symbolic dynamical system} associated to a language $L$ is the one-sided shift $S(x_0x_1x_2...)=x_1x_2...$ on the subset $X_L$ of ${\mathcal A}^{\GN}$ made with the infinite sequences such that for every $s'<s$, $x_{s'}...x_s$ is in $L$.

\noindent For a word $w=w_1...w_s$ in $L$, the {\em cylinder} $[w]$ is the set $\{x\in X_L; x_0=w_1, ... x_{s-1}=w_s\}$. \\
For a system $(X,T)$ and a finite partition $Z=\{Z_1,\ldots Z_{\rho}\}$ of $X$, the {\em trajectory} of a point $x$ in $X$ is the infinite
sequence
$(x_{n})\ind$ defined by $x_{n}=i$ if ${T}^nx$ falls into
$Z_i$, $1\leq i\leq \rho$.

\noindent Then $L(Z,T)$ is the language made of all the finite factors of all the  trajectories, and $X_{L(Z,T)}$ is the {\em coding} of $X$ by $Z$. 
\end{definition}

Note that the symbolic dynamical system $(X_L,S)$ is minimal (in the usual sense, every orbit is dense) if and only if the language $L$ is mimimal as in Definittion \ref{dln}.

     \begin{definition}\label{dlr}
A language $L$ or the symbolic system $(X_L,S)$ is {\em linearly recurrent}  if there exists $K$ such that in $L$, every word of length $n$ occurs in every word of length $Kn$. \end{definition}

\subsection{Boshernitzan's criteria for symbolic systems}
\begin{definition} For an invariant measure $\mu$ on $(X_L,S)$, 
let $e_n(S,\mu)$ be the smallest positive measure  of the cylinders  of length $n$. \end{definition}

The following sufficient condition for unique ergodicity is known as Boshernitzan's criterion; it is defined, named, and its sufficiency is proved for codings of interval exchange transformations  in \cite{vb}, then this is extended to every symbolic dynamical system in \cite{bo2}. 

\begin{pp}\label{boue} If $(X_L,S)$ is minimal,   the system is uniquely ergodic if  there exists an invariant measure  such that $\limsup_{n\to +\infty}ne_n(S,\mu)>0$. \end{pp}

The following result on linear recurrence is also due to M. Boshernitzan,  but was  written by T. Monteil in \cite{fogg}, Exercise 7.14.

\begin{pp}\label{bos} $(X_L,S)$ is linearly recurrent if and only if there exists an invariant measure on $(X_L,S)$ such that $\liminf_{n\to +\infty}ne_n(S,\mu)>0$. \end{pp}

\subsection{Measure-theoretic properties}
Let $(X, T, \mu )$ be a probability-preserving dynamical system.

\begin{definition}\label{drig} $(X, T, \mu )$ is {\em rigid} if there exists a sequence
    $q_{n}\to\infty$  such that for
any measurable set $A$ 
 $\mu (T^{q_{n}}A\Delta A)\to 0.$ \end{definition} 

\begin{definition}\label{nm} In $(X,T)$, a Rokhlin {\em tower}  is a collection 
    of disjoint measurable sets called {\em levels} $F$, $TF$, \ldots, $T^{h-1}F$.  $F$ is the {\em basis} of the tower.\\
    If $X$ is 
    equipped with a partition $P$ such that each level $T^{r}F$ is contained in
    one atom
$P_{w(r)}$, 
 the {\em name} of the 
tower is the word
 $w(0)\ldots w(h-1)$. \\
 A symbolic systems is generated by families of Rokhlin towers $F_{i,n}$, ..., $T^{h_{i,n}-1}F_{i,n}$, $1\leq i\leq K$, $n\geq 1$,  if each level in each towers is contained in a single atom of the partition into cylinders $\{x_0=i\}$, and for any word $W$ in $L(T)$ there exist $i$ and $n$ such that $W$ occurs in the name (for this partition) of the tower of basis $F_{i,n}$.
 \end{definition}
 
If a symbolic system is generated by families of Rokhlin towers, then, for any invariant measure, any measurable set can be approximated in measure by finite unions of levels of towers.

 \begin{definition}\label{dr1}  $(X, T,\mu )$  is of {\em rank one} if there exists a sequence of Rokhlin towers such that the whole $\sigma$-algebra is generated by the partitions $\{F_n, TF_n, \ldots, T^{h_n-1}F_n, X\setminus\cup_{j=0}^{h_n-1}T^jF_n\}$.
      \end{definition}
      
 \subsection{Rotations}
 The dynamical behavior of a rotation $R$ of angle $\alpha$ on the $1$-torus is linked with the {\em  Euclid continued fraction expansion of $\alpha$}. We assume the reader  is familiar with the notation  $\alpha=[0,a_1, a_2,...]$; we define in the classical way the convergents $\frac{p_n}{q_n}$ by $p_{-1}=1$, $q_{-1}=0$, $p_0=0$, $q_0=1$, $p_{n+1}=a_{n+1}p_n+p_{n-1}$, $q_{n+1}=a_{n+1}q_n+q_{n-1}$. Let $\alpha_n=\ab q_n\alpha -p_n\ab$. 
We recall

\begin{definition}
$ \alpha$ has {\em bounded partial quotients} if the $a_i$ are bounded.
 \end{definition}       

 {\em Throughout the paper,   except when we need more precision, we use $C$ as a generic notation for  constants}.

\section{Veech and Sataev  examples}

\begin{definition}\label{dv}
The {\em Veech 1969} system is   defined, if $Rx=x+\alpha$ modulo $1$, by $T(x,s)=(Rx,\sigma(x)s)$, $s=1,2$, where
\begin{itemize}
\item $\sigma(x)=\sigma_0$  if $x$ is in the interval $[0,\beta[$, $\sigma_0$ being the  exchange $E$,
\item $\sigma(x)=\sigma_1$  if $x$ is in the interval $[\beta,1[$, $\sigma_1$ being the identity $I$.
 \end{itemize}
 \end{definition}

This is defined  (in a slightly different terminology) in the famous paper \cite{ve69}, where Veech considers also the variant where  $\sigma_0=I$, $\sigma_1=E$. We can identify $[0,1[\times \{s\}$ with 
$[s-1,s[$; then $T$ is also an {\em interval exchange transformation } as in Figure 1 (note that six intervals appear in the picture, but two of them move together thus $T$ is indeed a $5$-interval exchange transformation).\\

\begin{figure}
\begin{center}
\begin{tikzpicture}[scale = 5]

\draw (0,0.55)node[above]{$0$};
\draw (.4,0.55)node[above]{$\beta$};
\draw (.75,0.55)node[above]{$1-\alpha$};
  \draw (1,0.55)node[above]{$1$};
\draw (1.1,0.55)node[above]{$1$};
\draw (1.5,0.55)node[above]{$1+\beta$};
\draw (1.85,0.55)node[above]{$2-\alpha$};
\draw (2.1,0.55)node[above]{$2$};

\draw(.4,.53)--(.4,.57);
\draw(.75,.53)--(.75,.57);
\draw(1.5,.53)--(1.5,.57);
\draw(1.85,.53)--(1.85,.57);
  
\draw(0,.55)--(.4,.55);
\draw(.4,.55)--(.75,.55);
\draw(.75,.55)--(1,.55);
\draw(1.1,.55)--(1.5,.55);
\draw(1.5,.55)--(1.85,.55);
\draw(1.85,.55)--(2.1,.55);

\draw (.2,0.55)node[below]{$1_1$};
\draw (.575,0.55)node[below]{$1_2$};
\draw (.875,0.55)node[below]{$1_3$};
  \draw (1.3,.55)node[below]{$2_1$};
\draw (1.675,0.55)node[below]{$2_2$};
\draw (1.975,0.55)node[below]{$2_3$};

\draw (0,0)node[above]{$0$};
\draw (.25,0)node[above]{$\alpha$};
\draw (.65,0)node[above]{$\alpha+\beta$};
  \draw (1,0)node[above]{$1$};
\draw (1.1,0)node[above]{$1$};
\draw (1.35,0)node[above]{$1+\alpha$};
\draw (1.75,0)node[above]{$1+\alpha+\beta$};
\draw (2.1,0)node[above]{$2$};

\draw(0,0)--(.25,0);
\draw(.25,0)--(.65,0);
\draw(.65,0)--(1,0);
\draw(1.1,0)--(1.35,0);
\draw(1.35,0)--(1.75,0);
\draw(1.75,0)--(2.1,0);

\draw(.25,-.02)--(.25,.02);
\draw(.65,-.02)--(.65,.02);
\draw(1.35,-.02)--(1.35,.02);
\draw(1.75,-.02)--(1.75,.02);

\draw (.125,0)node[below]{$T1_3$};
\draw (.45,0)node[below]{$T2_1$};
\draw (.825,0)node[below]{$T1_2$};
  \draw (1.225,0)node[below]{$T2_3$};
\draw (1.55,0)node[below]{$T1_1$};
\draw (1.925,0)node[below]{$T2_2$};

\end{tikzpicture}
\caption{Veech 1969}
\end{center}
\end{figure}
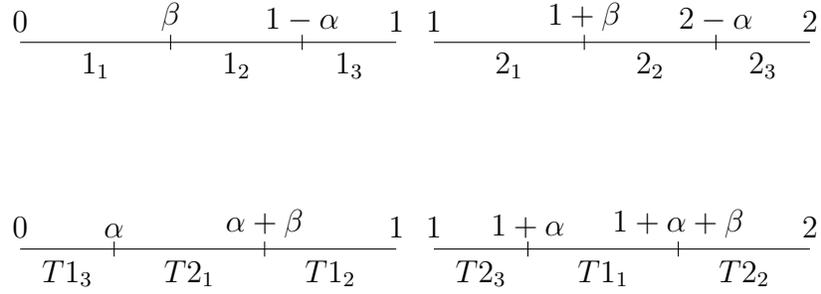

We can generalize Veech 1969 naturally by marking several points $\beta_i$, and taking more than two copies of the intervals: thus we take $r+1$ different permutations on $\{1,...,d\}$, changing permutation each time we cross a point $\beta_i$: these transformations are defined by Sataev \cite{sat} in 1975 for $d=r+1$, clearly without knowledge of Veech's work. Then these systems appear in \cite{ste} for $r=1$ and all $d$, and in \cite{mer} for all $r$ and $d$, but only in the particular case of  extensions by the commutative group $\GZ/d\GZ$, where all the permutations are circular; some of these systems are also considered in \cite{gp}.

 In general the $\beta_i$ will be chosen to be rationally independent from $\alpha$, but if there is only one $\beta$ and it is equal to $1-\alpha$, we get the {\em square-tiled interval exchange transformations} of \cite{fh2} (though the geometrical model is not the same); thus, to generalize both Veech 1969, Sataev 1975,   and the square-tiled interval exchange transformations, we keep the possibility of choosing one of the $\beta_i$ to be $1-\alpha$.   \\

Throughout this paper we take $\alpha$ irrational, $0<\beta_1<....<\beta_r<1$ irrational, with possibly $\beta_t=1-\alpha$; more precisely, {\em if the index $t$ exists, then 
$\beta_t=1-\alpha$; otherwise $\beta_i\neq 1-\alpha$ for all $i$}. We choose 
$\sigma_0$, ..., $\sigma_{r}$, permutations of $\{1,...,d\}$. We always suppose $\sigma_j\neq \sigma_{j+1}$, $0\leq j \leq r-1$,  as otherwise we could delete some $\beta_i$. We take all the $\beta_i$, $i\neq t$, and  all the $\beta_i-\beta_j$ not in $\mathbb Z(\alpha)$. We shall need sometimes another inequality,  which generalizes the non-commutation condition used in \cite{fh2}, and which we call the {\em product inequality}: namely, when $\beta_t=1-\alpha$ we ask that $\sigma_r\sigma_{t-1} \neq \sigma_0\sigma_t$, when $\beta_j\neq 1-\alpha$ for all $j$  we ask that $\sigma_r\neq \sigma_0$.

\begin{definition}\label{dvg}
The {\em generalized Veech - Sataev} system is defined, if $Rx=x+\alpha$ modulo $1$, by $T(x,s)=(Rx,\sigma(x)s)$, $1\leq s\leq d$, where 
\begin{itemize}
 \item $\sigma(x)=\sigma_j$ if $\beta_j \leq x <\beta_{j+1}$, $1\leq j\leq r-1$, 
\item $\sigma(x)=\sigma_0$ if $0 \leq x < \beta_1$, 
\item $\sigma(x)=\sigma_{r}$ if $\beta_{r} \leq x < 1$.\end{itemize}
\end{definition}

$T$ can be seen also as an interval exchange transformation on at most $d(r+1)$  intervals, or with the following geometric model, generalizing the {\em Masur-Smillie geometrical model} for Veech 1969 \cite{mat}: we build a surface by gluing $d$ tori, the interval $[\beta_i,\beta_{i+1}[$ in the right edge of the $s$-th torus being glued with the same interval in the left edge of the $\sigma_is$-th torus, and mutatis mutandis for the intervals $[0,\beta_{1}[$ and $[\beta_{r}, 1[$. Then   we take the directional flow of slope $\alpha$, going from one torus to the other when crossing the gluing lines, and  $T$ is its first return map on the union of the $d$ left vertical sides.

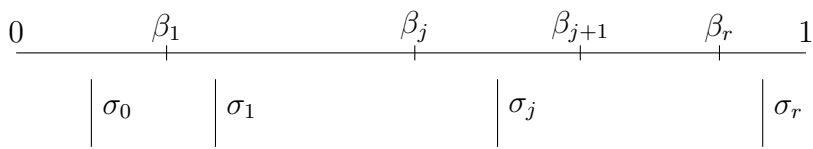
\begin{figure}
\begin{center}
\begin{tikzpicture}[scale = 5]

\draw (0,0.55)node[above]{$0$};
\draw (.4,0.55)node[above]{$\beta_1$};
\draw (1.06,0.55)node[above]{$\beta_j$};

\draw (1.5,0.55)node[above]{$\beta_{j+1}$};
\draw (1.87,0.55)node[above]{$\beta_r$};
\draw (2.1,0.55)node[above]{$1$};

\draw(.4,.53)--(.4,.57);
\draw(1.06,.53)--(1.06,.57);

\draw(1.5,.53)--(1.5,.57);
\draw(1.87,.53)--(1.87,.57);

\draw(0,.55)--(2.1,.55);
\draw(.2,.30)--(.2,.48);
\draw(.53,.30)--(.53,.48);
\draw(1.28,.30)--(1.28,.48);

\draw(1.985,.30)--(1.985,.48);

\draw (.2,.4)node[below, right]{$\sigma_0$};
\draw (.53,.4)node[below, right]{$\sigma_1$};
\draw (1.28,.4)node[below, right]{$\sigma_j$};
\draw (1.985,.4)node[below, right]{$\sigma_r$};

 \end{tikzpicture}
 \caption{Generalized Veech - Sataev}
\end{center}
\end{figure}

This model allows us to give a minimality condition for  generalized Veech - Sataev systems. 

\begin{pp}\label{pmn} If $\alpha$ and all the $\beta_i$ are irrational, all the $\beta_i$, $i\neq t$,  and all the $\beta_i-\beta_j$ are not in $\mathbb Z(\alpha)$, an NSC for minimality is that no strict subset of $\{1\ldots d\}$ is invariant by all the $\sigma_i$. \end{pp}
{\bf Proof}\\
If a strict subset $A$ of $\{1\ldots d\}$ is invariant by all the $\sigma_i$, then $\cup_{i\in A}[0,1[\times\{i\}$ is invariant par $T$, and $T$ is not minimal.

In the other direction, the condition on the permutations ensures that the surface defined above is connected, and the flow is minimal as the conditions on the $\beta_i$ ensure  there is no connection, except possibly (if $\beta_t=1-\alpha$) $d$ connections between $1-\alpha$ and $0$, each one staying inside one torus;  these connections do not separate the surface  into several parts.
\qed\\

\section{The rotation with marked points and the Ostrowski expansion}

\subsection{Rokhlin towers}\label{roto}

The rotation $R$ can be coded either by the partition $Z$ of the interval into $[0,1-\alpha[$ and $[1-\alpha, 1[$, or by the partition $Z'$ of the interval by the points $\beta_1,..,.\beta_r$. This gives two languages $L$ and $L'$, and two symbolic systems. The first one is the {\em natural coding} of $R$: it is assimilated to $R$ itself and denoted by $(X,R)$. The second one is called the {\em  rotation with marked points} and denoted by $(X',S)$.\\

\begin{figure}
\begin{center}
\begin{tikzpicture}[scale = 1]

\draw (1.5,3.5) node[above]{$(0)$};
\draw (1.5,9.5) node[above]{$(0)$};
\draw (4.5,3.5) node[above]{$(1)$};
\draw (10.5,3.5) node[above]{$(a_{n+1}-1)$};
\draw (12.5,3.5) node[above]{$(a_{n+1})$};

\draw(0,0)--(0,8);\draw(0,0)--(13,0);\draw(13,0)--(13,8);\draw(0,8)--(13,8);
\draw(0,8.6)--(0,11);\draw(0,11)--(3,11);\draw(3,8.6)--(3,11);\draw(0,8.6)--(3,8.6);
\draw[dashed](3,0)--(3,8);
\draw[dashed](6,0)--(6,8);\draw[dashed](9,0)--(9,8);\draw[dashed](12,0)--(12,8);

\draw(10,0) node[above]{$\alpha$};
\draw(10,-0.8) node[above]{$0$};

\draw(10,-0.1)--(10,0.1);
\draw(10,-0.7)--(10,-0.9);

\draw(0,11) node[left]{$0$};
\draw(3,11) node[above]{$\alpha_n$};
\draw(0,8) node[left]{$-\alpha_{n-1}$};
\draw(0,8.6) node[left]{$\alpha-\alpha_{n-1}$};
\draw(13,-0.8) node[right]{$\alpha_n$};

\draw(0,0) node[below]{$\alpha-\alpha_{n-1}+\alpha_n$};

\draw(3,7.8) node[right]{$-\alpha_{n-1}+\alpha_n$};

\draw(13,8) node[right]{$0$};
\draw(12,8) node[above]{$-\alpha_{n+1}$};
\draw[dashed](0,-0.8)--(13,-0.8);

\draw(13,0) node[right]{$\alpha_n+\alpha$};
\draw[dashed](13,-0.8)--(13,-1.8);

\draw[dashed](10,-0.8)--(10,-1.8);

\draw(0,10.4) node[left]{$1-\alpha$};
\draw(13,7.4) node[right]{$1-\alpha$};
\draw(-0.1,10.4)--(0.1,10.4);
\draw(12.9,7.4)--(13.1,7.4);

\end{tikzpicture}

\caption{Rokhlin $n$-towers for the rotation, $n$ odd}
\end{center}
\end{figure}
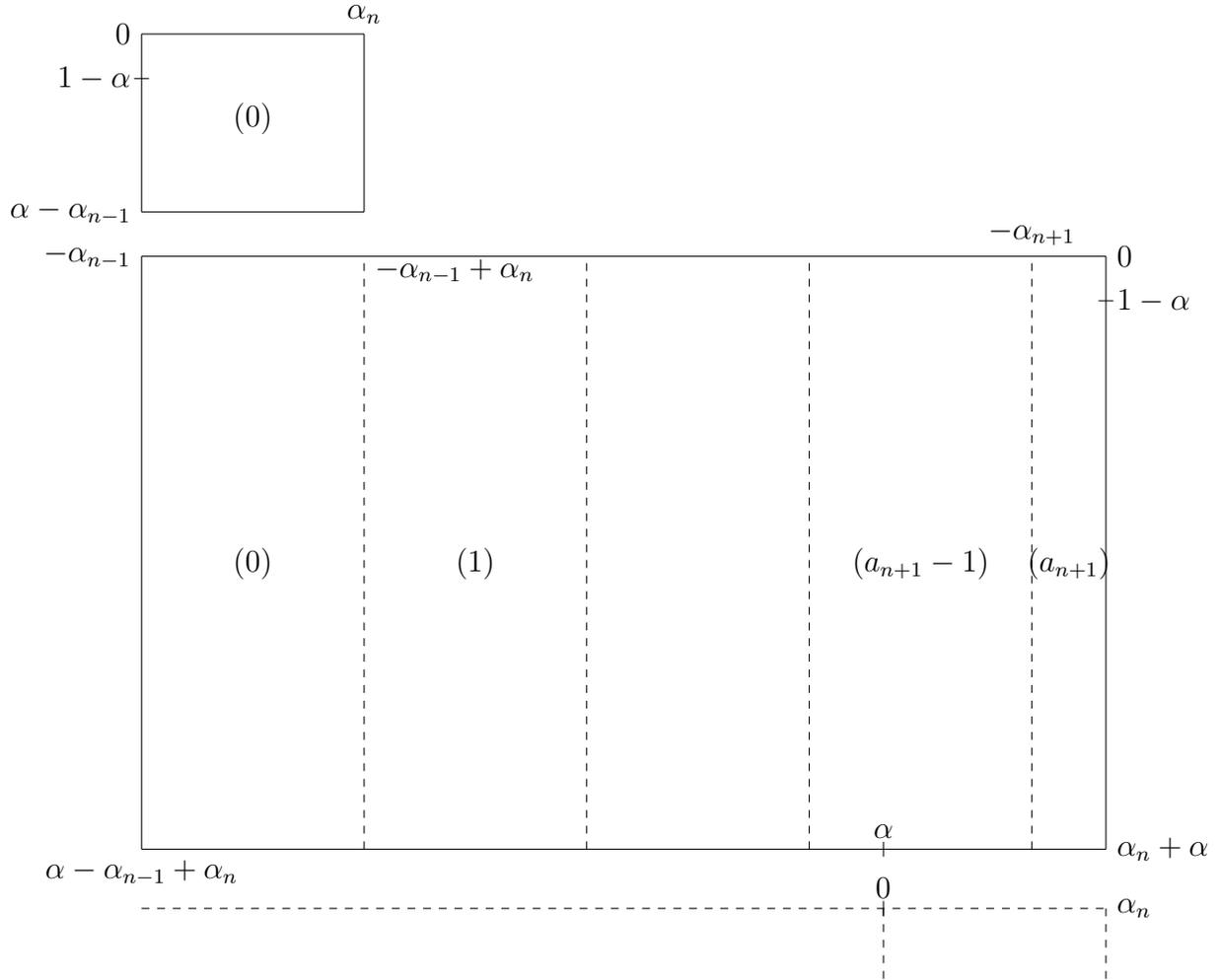

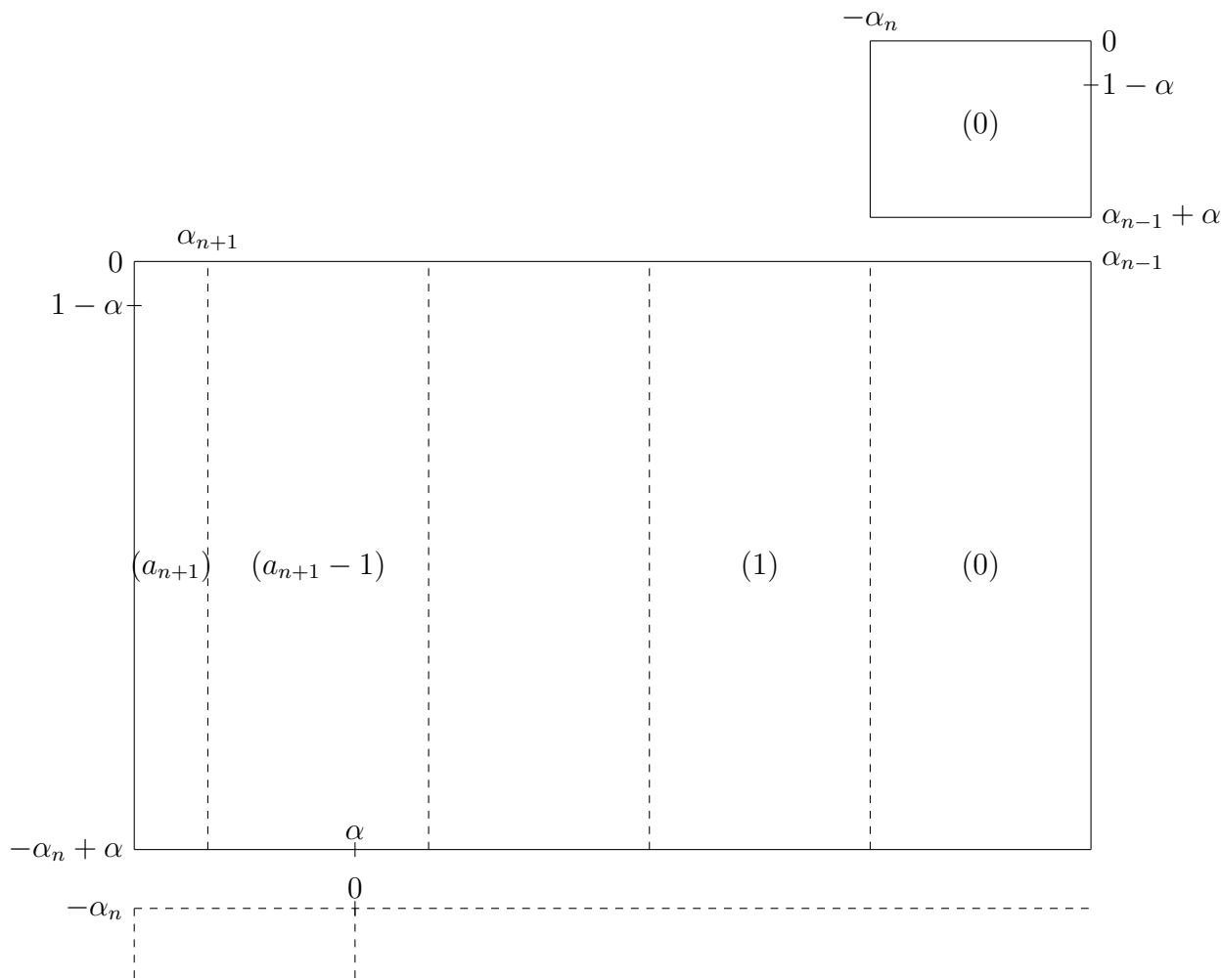
\begin{figure}
\begin{center}
\begin{tikzpicture}[scale = 1]

\draw (11.5,3.5) node[above]{$(0)$};
\draw (11.5,9.5) node[above]{$(0)$};
\draw (8.5,3.5) node[above]{$(1)$};
\draw (2.5,3.5) node[above]{$(a_{n+1}-1)$};
\draw (.5,3.5) node[above]{$(a_{n+1})$};

\draw(13,0)--(13,8);\draw(0,0)--(13,0);\draw(0,0)--(0,8);\draw(0,8)--(13,8);
\draw(13,8.6)--(13,11);\draw(10,11)--(13,11);\draw(10,8.6)--(10,11);\draw(10,8.6)--(13,8.6);
\draw[dashed](10,0)--(10,8);
\draw[dashed](7,0)--(7,8);\draw[dashed](4,0)--(4,8);\draw[dashed](1,0)--(1,8);

\draw(3,0) node[above]{$\alpha$};
\draw(3,-0.8) node[above]{$0$};

\draw(3,-0.1)--(3,0.1);
\draw(3,-0.7)--(3,-0.9);

\draw(13,11) node[right]{$0$};
\draw(10,11) node[above]{$-\alpha_n$};
\draw(13,8) node[right]{$\alpha_{n-1}$};
\draw(13,8.6) node[right]{$\alpha_{n-1}+\alpha$};
\draw(0,-0.8) node[left]{$-\alpha_n$};

\draw(0,8) node[left]{$0$};
\draw(1,8) node[above]{$\alpha_{n+1}$};
\draw[dashed](0,-0.8)--(13,-0.8);

\draw(0,0) node[left]{$-\alpha_n+\alpha$};
\draw[dashed](3,-0.8)--(3,-1.8);

\draw[dashed](0,-0.8)--(0,-1.8);

\draw(0,7.4) node[left]{$1-\alpha$};
\draw(13,10.4) node[right]{$1-\alpha$};
\draw(-0.1,7.4)--(0.1,7.4);
\draw(12.9,10.4)--(13.1,10.4);

\end{tikzpicture}

\caption{Rokhlin $n$-towers for the rotation, $n$ even}
\end{center}
\end{figure}

It is well known, and written for example in \cite{gp}, that for the rotation $R$ its natural coding is generated by two families of Rokhlin towers, made of intervals. We shall now describe precisely the towers at stage $n$, or $n$-towers.

At each stage $n\geq 1$, there are  one {\em large tower} made of $q_n$ intervals (or {\em levels}) of length $\alpha_{n-1}$ and one {\em small tower} made of $q_{n-1}$ intervals (or  levels) of length $\alpha_{n.}$. These are described in 
 Figure 3 if $n$ is odd, we make all our comments in that case; the case when  $n\geq 2$ is even  can be  deduced, mutatis mutandis, from Figure 4. Namely, the  large tower is represented by the lower rectangle, and the small tower by the upper rectangle.The rotation $R$ sends the basis $[-\alpha_{n-1}+\alpha_n+\alpha,\alpha_n+\alpha[$ to an interval which we put just above it, and call a level of the large tower; this interval is sent by $R$ just above, and so on until, by $R^{q_n-1}$ applied to the basis of the large tower,  we reach the top of the large tower, $[-\alpha_{n-1},0[$. Then the left part of this top, $[-\alpha_{n-1},-\alpha_{n-1}+\alpha_n[$ is sent by $R$ onto he basis 
 $[-\alpha_{n-1}+\alpha,-\alpha_{n-1}+\alpha_n+\alpha[ $ of the small tower, and we go up in the small tower until, by $R^{q_{n-1}-1}$ 
 applied to the basis of the small tower, we reach the top of the small tower, $[0,\alpha_n[$. Where we go next from $[0,\alpha_n[$ or $[-\alpha_{n-1}+\alpha_n,\alpha_n[$ is shown at the bottom of the picture, one application of $R$ goes to the point just above, in the basis of the large tower. For any $x$, $R^{q_n}x$ is the point situated at distance $\alpha_n$ to the left of $x$ (extending the intervals if one of these points is not in the picture). Note that  the three points $1-\alpha$, $0$ and $\alpha$ can be considered as very close together in the $n$-towers. 
 
 Each level of each tower is included in one atom of the partition $Z$.  At the  beginning, if $\alpha >\frac12$, the large $1$-tower has one level, the interval $[1-\alpha, 1[$ and he small $1$-tower has one level, the interval $[0,1-\alpha[$, Figure 3 is still valid. If $\alpha <\frac12$, the $1$-towers, which are still given by Figure 3, are more complicated, but we can define $0$-towers:
 the large $0$-tower has one level, the interval $[0, 1-\alpha[$ and he small $0$-tower has one level, the interval $[1-\alpha,1[$.

The large tower is partitioned from left to right into $a_{n+1}+1$ columns, of width $\alpha_n$ except for the last one which is of width $\alpha_{n+1}$, and except for $\alpha <\frac12$, $n=0$, where there are only $a_1$ columns and thus Figure 4 does not apply. We denote the columns as in Figure 3 or 4, and include the whole small tower in column $(0)$. The description of $R$ defines immediately the next towers: to get the large $n+1$-tower we stack the columns  of the $n$-tower above each other, with the column $(a_{n+1}-1)$ at the  bottom, then $(a_{n+1}-2)$, ..., $(0)$, and the small $n$-tower at the top, while the $n$-column $(a_{n+1})$ becomes the small $n+1$-tower.\\

{\em Note that all levels are semi-open intervals, closed on the left, open on the right, and thus each column includes its left vertical side and not its right one}.

The following lemma will be fundamental in our computations: we shall use it when $\alpha$ has bounded partial quotients, as it quantifies the linear recurrence of the natural coding of the rotation,  but we can state and prove it in the general case. 

\begin{lemma}\label{traj} Suppose $x$ or $y$, or both, are in the basis of the large $n$-tower, and $x=y+z$, $\alpha_{n+1}\leq z\leq \alpha_n$. Then the smallest $k>0$ such that $\alpha$ lies between $R^kx$ and $R^ky$, with $\alpha\neq R^kx$,   is at least $q_n+q_{n-1}$ and at most $q_{n+2}+q_{n+1}+q_n$. \end{lemma}
{\bf Proof}\\
$y$ is at  a distance $z$  to the left of $x$; then for all $m$ $R^my$ is at the same distance of $R^mx$ on the circle.
We make all computations with $n$ odd, the even case is similar. To simplify notations, we write them for $x$ and $y$ which are not on the sides of any $n$-columns, thus excluding a countable set. However, we  notice they are still valid on this countable set,  because of our conventions that the columns are closed on the left, open on the right and we allow $\alpha= R^ky$ when saying $\alpha$ appears between the two orbits.

\begin{itemize}
\item{(i)} Suppose first $x$ is in the basis of the large $n$-tower and  not $y$. Then $x$ is at a distance at most $\alpha_n$ from the left of the large tower, thus in column  $(0)$; $y$ is to the left of the large tower and at less than $\alpha_n$ from it, thus between $-\alpha_{n-1}+\alpha$ and $-\alpha_{n-1}+\alpha_n+\alpha$, thus in the basis of the small tower (see Figure 9 below). $y$ is at a distance $d_1$ from the left of this basis, $x$ is at a distance $0<d_2<z$ from the left of the large tower, with $d_1+z=d_2+\alpha_n$. We make $q_n+q_{n-1}$ iterations of $R$. The orbit of $x$ goes up through the large and small towers,  and at the  $q_n+q_{n-1}$-th iteration hits the basis of the large tower, at a point situated $d_2$ to the right of $\alpha$, The orbit of $y$ goes up through the small tower, at the $q_{n-1}$-th iteration hits the basis of the large tower at a point situated  $d_1$ to the right of $\alpha$, then at the $q_n+q_{n-1}$-th iteration hits this basis again, $\alpha_n$ left of the previous hit, thus left of $\alpha$ as $d_1-\alpha_n=d_2-z<0$; and before the  $q_n+q_{n-1}$-th iteration $\alpha$ does not appear between the two orbits.
\item{(ii)} Suppose $x$ and $y$ are in the basis of the large tower and in two different columns. Then these columns must be  adjacent, and, after at most $a_{n+1}q_n$ iterations of $R$, during which $\alpha$ does not appear between the two orbits, we are in the situation of   case $(i)$. 
\item{(iii)} Suppose $x$ and $y$ are in the basis of the large tower and in the same column. This column cannot be column $(a_{n+1})$, and $x$ is at distance $d_3$ from the right of its column. After at most $a_{n+1}q_n+q_{n-1}=q_{n+1}$ iterations of $R$, during which $\alpha$ does not appear between the two orbits, the orbit of $x$ hits the basis of the large tower, at a point situated $d_3$ from its right end, and $y$ also, at distance $d_3+z$ from the right. At this moment, if $d_3<\alpha_{n+1}$, then the orbit of $x$ is in column $(a_{n+1})$ and  the orbit of $y$ is in column $(a_{n+1}-1)$ and we are in the situation of case $(ii)$. If $d_3>\alpha_{n+1}$,   the orbits of $x$ and $y$ are in column $(a_{n+1}-1)$, and we are again in case $(iii)$, but with $d_3$ replaced by $d_3-\alpha_{n+1}$. As $d_3<\alpha_n$ and $\alpha_n=a_{n+2}\alpha_{n+1}+\alpha_{n+2}$, after at most $a_{n+2}$ such laps,  during which $\alpha$ does not appear between the two orbits, we are in the situation of   case $(ii)$. 
\item{(iv)} Suppose finally  $y$ is in the basis of the large $n$-tower and  not $x$. Then $x$ is to the right of the large tower and at less than $\alpha_n$ from it, and $y$ is to the right of $\alpha$; after  $q_n$ iterations of $R$, during which $\alpha$ does not appear between the two orbits, we are in the situation of   case $(ii)$ or $(iii)$. 
\end{itemize}

Thus, by taking  case $(i)$ for the minimum, and summing our estimates for the maximum, we get the required result.\qed\\

\begin{corollary}\label{ctraj} Let $\beta$ and $\beta'$ be any two points on the circle. Suppose  $x=y+z$, $\alpha_{n+1}\leq z\leq \alpha_n$, and $\beta$ lies between $x$ and $y$, with $\beta\neq x$. Then the smallest $k>0$ such that $\beta'$ lies between $R^kx$ and $R^ky$, with $\beta'\neq R^kx$,   is at most $q_{n+2}+q_{n+1}+q_n$. \end{corollary}
{\bf Proof}\\ By Lemma \ref{traj} this is true for $\beta'=\alpha$ and  any $\beta$, just because $\beta$ is in one of the $n$-towers and is the image of some point in the basis of the large one. As $R$ commmutes with every translation, this is true also for any $\beta$ and $\beta'$. \qed\\

\subsection{Ostrowski expansion}
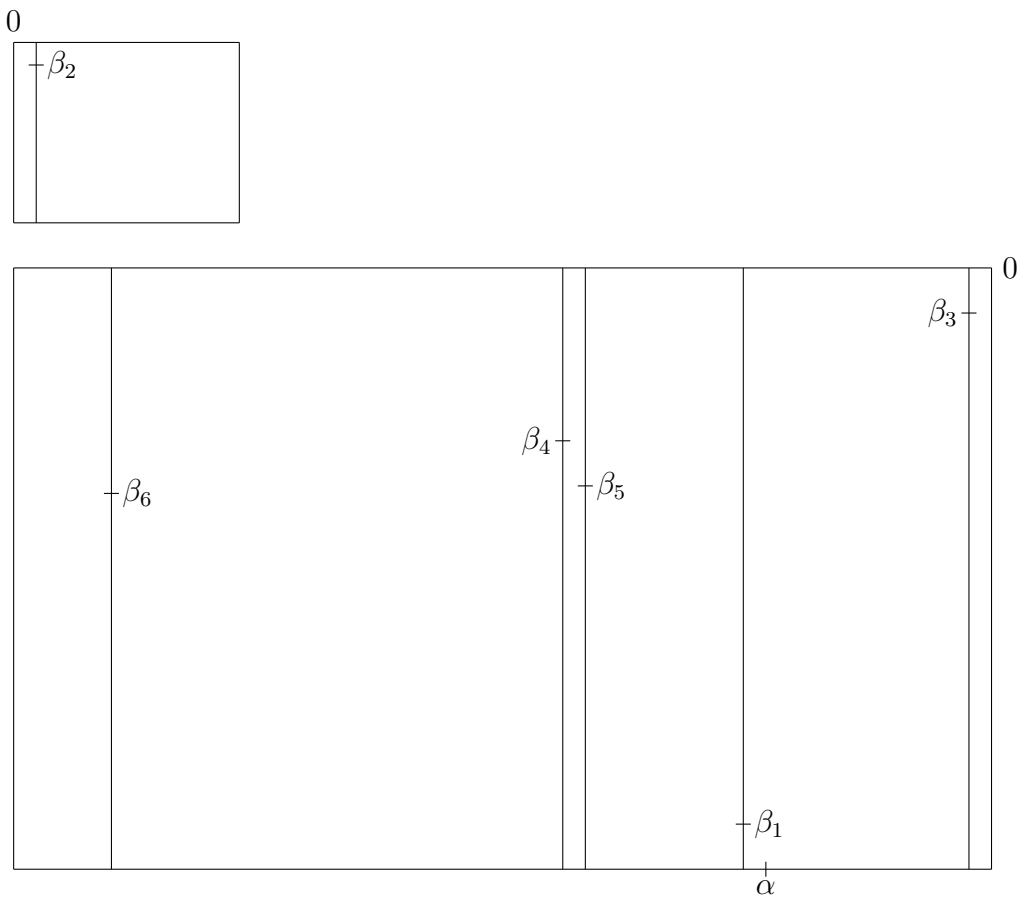
\begin{figure}
\begin{center}

\begin{tikzpicture}[scale = 1]

\draw(0,0)--(0,8);\draw(0,0)--(13,0);\draw(13,0)--(13,8);\draw(0,8)--(13,8);
\draw(0,8.6)--(0,11);\draw(0,11)--(3,11);\draw(3,8.6)--(3,11);\draw(0,8.6)--(3,8.6);

\draw(10,0) node[below]{$\alpha$};
\draw(9.7,0.6) node[above, right]{$\beta_1$};
\draw(0.3,10.7) node[above, right]{$\beta_2$};
\draw(12.7,7.4) node[below, left]{$\beta_3$};
\draw(7.3,5.7) node[above, left]{$\beta_4$};
\draw(7.6,5.1) node[above, right]{$\beta_5$};
\draw(1.3,5) node[above, right]{$\beta_6$};

\draw(10,-0.1)--(10,0.1);
\draw(9.6,0.6)--(9.8,0.6);
\draw(0.2,10.7)--(0.4,10.7);
\draw(12.6,7.4)--(12.8,7.4);
\draw(7.2,5.7)--(7.4,5.7);
\draw(7.5,5.1)--(7.7,5.1);
\draw(1.2,5)--(1.4,5);

\draw(1.3,0)--(1.3,8);
\draw(0.3,8.6)--(0.3,11);
\draw(7.3,0)--(7.3,8);
\draw(7.6,0)--(7.6,8);
\draw(9.7,0)--(9.7,8);
\draw(12.7,0)--(12.7,8);

\draw(0,11) node[above]{$0$};

\draw(13,8) node[right]{$0$};

\end{tikzpicture}

\caption{Rokhlin $n$-towers for the rotation with marked points}
\end{center}
\end{figure}

We now put the points $\beta_i$, $i\neq t$,  in the picture. By partitioning the two towers for $R$ as in Figure 5  (for odd $n$), we get $r+2$ towers generating the rotation with marked points $S$, for which each level of each tower is included in one atom of the partition $Z'$ (if $\beta_t=1-\alpha$, only $r+1$ towers are needed as $1-\alpha$ is on the side of one tower).

For each $1\leq i\leq r$, $i\neq t$, and $n\geq 1$ we define $b_{n+1}(\beta_i)$ as an integer between $0$ and $a_{n+1}$.

\begin{definition}\label{dost} $b_{n+1}(\beta_i)$
 is $b\neq 0$ if $\beta_i$ is in column $(b)$ of the large $n$-tower (for $R$), and $0$ if $\beta_i$ is either in column $(0)$ of the large $n$-tower  or in the small $n$-tower.\\
 For odd $n$ (resp. even $n\geq 2$) let $x_{n}(\beta_i)$ be the (positive) distance of $\beta_i$ to the left (resp. right) side of the large and small $n$-towers in Figure 5 (resp. 6). \end{definition}

\begin{pp}\label{dostr}
For each $i\neq t$, the $b_n(\beta_i)$ are given by a form of  alternating  Ostrowski expansion of  $\beta_i$ by $\alpha$, where the Markovian condition is $b_n(\beta_i)=a_n$ implies $b_{n+1}(\beta_i)=0$. For $i\neq t$, $\beta_i$ is in $\mathbb Z(\alpha)$ if and only if either $b_n(\beta_i)=a_n-1$ for all $n$ large enough, or $b_{2n}(\beta_i)=a_{2n}$ for all $n$ large enough, or $b_{2n+1}(\beta_i)=a_{2n+1}$ for all $n$ large enough. For $i\neq t$, $j\neq t$, $\beta_i-\beta_j$ is in $\mathbb Z(\alpha)$ if and only if  $b_n(\beta_i)=b_n(\beta_j)$ for all $n$ large enough.\end{pp}
{\bf Proof}\\
We fix an $i\neq t$.  Then 
$$b_{n+1}(\beta_i)=\left[\frac{x_n(\beta_i)}{\alpha_n}\right].$$ 
Now, $x_{n+1}(\beta_i)$ is the distance of $\beta_i$ to the right (resp. left) side of its $n$-column if $n$ is odd (resp. even). Thus  we get  $x_{n}(\beta_i)=b_{n+1}(\beta_i)\alpha_n+\alpha_n-x_{n+1}(\beta_i)$ if $\beta_i$ is not in column $(a_{n+1})$, 
 $x_{n}(\beta_i)=b_{n+1}(\beta_i)\alpha_n+\alpha_{n-1}-x_{n+1}(\beta_i)$ if $\beta_i$ is in column $a_{n+1}$.  Note that if $\beta_i$ is in column $(a_{n+1})$ in the large $n$-tower, then it is in the small $n+1$-tower. Thus $b_{n+1}(\beta_i)=a_{n+1}$ implies $b_{n+2}(\beta_i)=0$, and this is the only Markovian condition they have to satisfy.
 
Thus when  $x_{n}(\beta_i)=b_{n+1}(\beta_i)\alpha_n+\alpha_{n-1}-x_{n+1}(\beta_i)$, then $x_{n+1}(\beta_i)=\alpha_{n+1}-x_{n+2}(\beta_i)$ and 
$x_{n}(\beta_i)=a_{n+1}\alpha_n-x_{n+2}(\beta_i)$. Together with the  formula when $b_{n+1}(\beta_i)<a_{n+1}$, this gives an expansion  $x_{1}(\beta_i)=\sum_{n\geq 1}(-1)^{n+1}\bar b_{n+1}\alpha_n$ with $\bar b_n=b_n(\beta_i)+1$ if $b_n(\beta_i)<a_{n}$, $\bar b_n=b_n$ if $b_n(\beta_i)=a_{n}$. Thus the $\bar b_n$ satisfy he Markovian condition $\bar b_{n-1}=a_{n-1}$ if $\bar{b_n}=0$.

Thus we identify the $\bar b_n$ with the alternating Ostrowski expansion of $x_1(\beta_i)$ by $\alpha$ defined in \cite{afh}.  If $\alpha>\frac12$,  $x_1(\beta_i)$ is either $\beta_i$ or $\beta_i+\alpha -1$; if $\alpha<\frac12$, using the $0$-towers, we can define $0\leq b_1(\beta_i) \leq a_1-1$ and $x_0(\beta_i)$ in the usual way, so that 
$x_{1}(\beta_i)=-x_0(\beta_i)+((b_1(\beta_i) +1)\alpha) \wedge (1-\alpha)$, and $x_0(\beta_i)$ is either $1-\beta_i$ or $1-\alpha - \beta_i$. In both cases, we get an expansion of $\beta_i$ by $\alpha$, which is $\beta_{i}=\sum_{n\geq 0}(-1)^{n+1}\bar b_{n+1}\alpha_n$ with a suitable $\bar b_1$, thus our   $b_n(\beta_i)$ do provide a form of  alternating Ostrowski expansion of $\beta_i$ by $\alpha$. 

The last conditions come from the fact that if $\beta_i=R^k\alpha$ for $k>0$ then $\beta_i$ is in the same column as $\alpha$, namely column $a_n-1$, in the $n-1$-towers for all $n$ large enough, while if 
$\beta_i=R^k\alpha$ for $k<0$ then $\beta_i$ is in the same column as $0$, and this alternates between $0$ (in the small tower) and $a_n$, and in both cases the converse is true by construction of the towers, as  the vertical distance from $\beta_i$ to $\alpha$ (resp. $0$) in the $n-1$-towers is ultimately constant while the horizontal distance tends to zero with $n$. Similarly, $\beta_i=R^k\beta_j$ if and only if in the $n-1$-towers $\beta_i$ is in the same column as $\beta_j$ for all $n$ large enough.
\qed\\

As a consequence, we can build $\beta_i$ with any prescribed sequence $0\leq b_n(\beta_i)\leq a_n$ satisfying the Markovian condition.

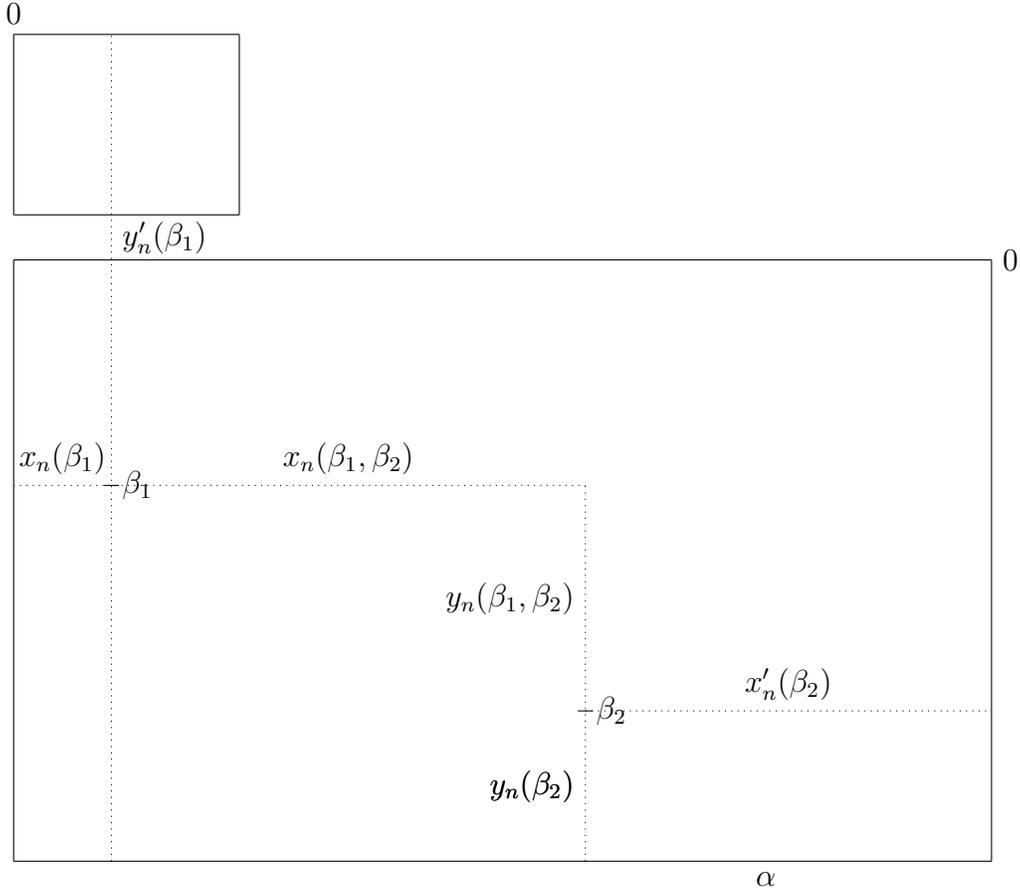
\begin{figure}
\begin{center}

\begin{tikzpicture}[scale = 1]

\draw(0,0)--(0,8);\draw(0,0)--(13,0);\draw(13,0)--(13,8);\draw(0,8)--(13,8);
\draw(0,8.6)--(0,11);\draw(0,11)--(3,11);\draw(3,8.6)--(3,11);\draw(0,8.6)--(3,8.6);

\draw(10,0) node[below]{$\alpha$};

\draw(7.6,2) node[above, right]{$\beta_2$};
\draw(1.3,5) node[above, right]{$\beta_1$};

\draw(7.5,2)--(7.7,2);
\draw(1.2,5)--(1.4,5);

\draw[dotted](0,5)--(7.6,5);
\draw[dotted](7.6,0)--(7.6,5);
\draw[dotted](1.3,0)--(1.3,11);
\draw[dotted](7.6,2)--(13,2);

\draw(0,11) node[above]{$0$};

\draw(13,8) node[right]{$0$};

\draw(7.6,1) node[left]{$y_n(\beta_2)$};
\draw(1.3,8.3) node[right]{$y'_n(\beta_1)$};
\draw(7.6,1) node[left]{$y_n(\beta_2)$};
\draw(7.6,3.5) node[left]{$y_n(\beta_1,\beta_2)$};

\draw(4.45,5) node[above]{$x_n(\beta_1,\beta_2)$};
\draw(.65,5) node[above]{$x_n(\beta_1)$};
\draw(10.3,2) node[above]{$x'_n(\beta_2)$};

\end{tikzpicture}

\caption{Positioning   the $\beta_i$ in the $n$-towers}
\end{center}
\end{figure}

\subsection{Linear recurrence}

To prove the  next theorem, we need some new notations. 

\begin{definition} For a given $n$,  each $\beta_i$, $i\neq t$, appears  in a single position in the $n$-towers as in Figure 6; it is   determined by $x_n(\beta_i)$, from Definition \ref{dost}. We shall use also
\begin{itemize}
\item $y_n(\beta_i)=y$ if $\beta_i=R^y\beta'_i$ where $\beta'_i$ is in the basis of the large $n$-tower, 
\item 
$x'_n(\beta_i)=\alpha_{n-1}-x_n(\beta_i)$, 
\item $x_n(\beta_i,\beta_j)=x_n(\beta_j,\beta_i)=\ab x_n(\beta_i)-x_n(\beta_j) \ab$,
\item  $x_n(\beta_i,\alpha)=x_n(\alpha,\beta_i)=\ab x'_n(\beta_i)-\alpha \ab$.,
\item  $y'_n(\beta_i)=q_n-y_n(\beta_i)$ if $\beta_i$ is in the large $n$-tower,  $y'_n(\beta_i)=q_n+q_{n-1}-y_n(\beta_i)$ if $\beta_i$ is in the small $n$-tower, 
\item when $\beta_i$ is in the small tower, $y"_n(\beta_i)=y_n(\beta_i)-q_n$,
\item $y_n(\beta_i,\beta_j)= y_n(\beta_j,\beta_i)=\ab y_n(\beta_i)-y_n(\beta_j) \ab$.
\end{itemize}\end{definition}

It is worth mentioning that $\beta_i$ and $\beta_j$ are close to each other in the $n$-towers vertically either if $y_n(\beta_i, \beta_j)$ is small or if $y_n(\beta_i) +y'_n(\beta_j)$ is small, and $\beta_i$ and $\beta_j$ are close to each other in the $n$-towers horizontally either if $x_n(\beta_i, \beta_j)$ is small or if $x_n(\beta_i) +x'_n(\beta_j)$ is small. Though this will not be mentioned explicitly, each time we claim $\beta_i$ and $\beta_j$ are far from each other in one of these senses, this means that we have checked both conditions.

\begin{theorem}\label{ostlr} The symbolic system $(X',S)$ is  linearly recurrent if and only all the following conditions are satisfied \begin{itemize}
\item $\alpha$ has  bounded partial quotients,
\item for each $i\neq t$, the number of consecutive $n$ such that $b_n(\beta_i)=a_n-1$ is bounded,
\item for each $i\neq t$, the number of consecutive $n$ such that $b_{2n}(\beta_i)=a_{2n}$ and the number of consecutive $n$ such that $b_{2n+1}(\beta_i)=a_{2n+1}$ are bounded,
\item for each  $i\neq t$, $j\neq t$ with $j\neq i$, the number of consecutive $n$ such that $b_n(\beta_j)=b_n(\beta_i)$ is bounded.
\end{itemize}
\end{theorem}
{\bf Proof}\\
We suppose first our conditions are not satisfied.

If $\alpha$ has unbounded partial quotients, there exists $n$ such that $q_n\alpha_n$ is arbitrarily small. In the  large $n$-tower, with $\alpha$ in the basis and $1-\alpha$ just below $0$, we see a cylinder, for the natural coding, of length $q_n-1$ and Lebesgue measure $\alpha_n$; this is a union of cylinders for the coding with marked points, of the same length and of smaller measure. As the Lebesgue measure is the only invariant measure by $R$, this contradicts linear recurrence by Proposition \ref{bos}.

 If $b_{n+1}(\beta_i)=a_{n+1}-1$, then by construction of the towers $x_{n+1}(\beta_i,\alpha) = x_n(\beta_i,\alpha)$ and $y_{n+1}(\beta_i)=y_n(\beta_i)$. If this holds for all  $M \leq n\leq M+N$, then $x_{M}(\beta_i,\alpha) = x_{M+N}(\beta_i,\alpha)\leq \alpha_{M+N}$ and $y_{M+N}(\beta_i)=y_M(\beta_i)\leq q_M$. Thus, for example, in the large $M+N$-tower we see a cylinder (for the coding with marked points) of measure $x_{M+N}(\beta_i,\alpha)$ and length $y_{M+N}(\beta_i)$. The product of these quantities  is at most $q_M\alpha_{M+N}\leq \theta^{-N}q_{M+N}\alpha_{M+N}\leq \theta^{-N}C$, where $\theta$ is the golden ratio. Thus, if $N$ is allowed to be arbitrarily large, this contradicts linear recurrence by Proposition \ref{bos}.
 
 If $b_{n+1}(\beta_i)=a_{n+1}$, then $\beta_i$ is in the small $n+1$-tower, with $x'_{n+1}(\beta_i) = x_n(\beta_i)$ and $y'_{n+1}(\beta_i)=y'_n(\beta_i)$, and then $x'_{n+1}(\beta_i) = x_{n+2}(\beta_i)$ and $y'_{n+1}(\beta_i)=y'_{n+2}(\beta_i)$. If this holds for all  $M \leq n\leq M+2N$, $y'_{M+2N}(\beta_i)=y'_M(\beta_i)\leq q_M$. In the large $M+2N$-tower we see a cylinder of measure $x_{M+2N}(\beta_i)\leq \alpha_{M+2N}$ and length $y'_{M+2N}(\beta_i)$. If $N$ is allowed to be arbitrarily large, we conclude as in the previous case.
 
  If $b_{n+1}(\beta_i)=b_{n+1}(\beta_j)$, then  $x_{n+1}(\beta_i,\beta_j) = x_n(\beta_i,\beta_j)$ and $y_{n+1}(\beta_i,\beta_j)=y_n(\beta_i,\beta_j)$. If this holds for all  $M \leq n\leq M+N$, then  $y_{M+N}(\beta_i,\beta_j)=y_M(\beta_i,\beta_j)\leq q_M$. In the  $M+N$-towers we see a cylinder of measure $x_{M+N}(\beta_i,\beta_j)\leq \alpha_{M+N}$ and length $y_{M+N}(\beta_i,\beta_j)$. If $N$ is allowed to be arbitrarily large, we conclude as in the previous cases.\\
  
  We suppose now all our conditions are satisfied. In particular, $\alpha$ has bounded partial quotients.
  
  If $b_{n+1}\beta_i)\neq a_{n+1}-1$, then $x_n(\beta_i,\alpha)\geq \alpha_{n+2}$. Otherwise, $x_n(\beta_i,\alpha)=x_m(\beta_i,\alpha)$ for the first $m>n$ for which $b_{m+1}(\beta_i)\neq a_{m+1}-1$, and we know $m\leq n+K$. Thus we get that for all $n$, 
   $$x_n(\beta_i,\alpha)\geq \alpha_{n+K+2}\geq C\alpha_n.$$
  If $b_{n}\beta_i)\neq a_{n}-1$, then by construction of the towers $y_n(\beta_i)\geq q_{n-1}$. Otherwise, $y_n(\beta_i)=y_m(\beta_i)$ for the last $m<n$ for which $b_{m}(\beta_i)\neq a_{m}-1$, and we know $m\leq n-K$. Thus we get that for all $n$, 
   $$y_n(\beta_i)\geq q_{n-K-1}\geq Cq_n.$$
   If $b_{n-1}\beta_i)\neq a_{n-1}-1$, by construction of the towers the result $y_{n-1}(\beta_i) \geq Cq_{n-1}$ implies, when $y"_n(\beta_i)$ is defined, that
   $$y"_n(\beta_i)\geq Cq_n.$$
    If $b_{n+1}\beta_i)\neq a_{n+1}$, then $x'_n(\beta_i)\geq \alpha_{n+1}$. Otherwise, $x'_n(\beta_i)=x'_m(\beta_i)$ for the first $m>n$ such that $m-n$ is even and  $b_{m+1}(\beta_i)\neq a_{m+1}$, and we know $m\leq n+K$. Thus we get that for all $n$, 
   $$x'_n(\beta_i)\geq \alpha_{n+K+1}\geq C\alpha_n.$$
If $b_{n+2}\beta_i)\neq a_{n+2}$, then $x'_{n+1}(\beta_i)\geq \alpha_{n+2}$ and by construction of the towers $x_n\geq\alpha_{n+2}$ ($\beta_i$ being far from one side of the $n+1$-towers, is far from the opposite side of the $n$-towers).  Otherwise, $x'_{n+1}(\beta_i)=x'_{m+1}(\beta_i)$ for the first $m>n$ such that $m-n$ is even and  $b_{m+2}(\beta_i)\neq a_{m+2}$, and we know $m\leq n+K$. Thus we get that for all $n$, 
   $$x_n(\beta_i)\geq \alpha_{n+K+2}\geq C\alpha_n.$$
       If $b_{n-1}\beta_i)\neq a_{n-1}$, then $\beta_i$ is not in the small $n-1$-tower, thus far from the top in the $n$-towers: we have $y'_n\geq q_{n-1}$. Otherwise, $y_n(\beta_i)=y_m(\beta_i)$ for the last $m<n$ such that $n-m$ is even and  $b_{m-1}(\beta_i)\neq a_{m-1}$, and we know $m\leq n-K$. Thus we get that for all $n$, 
   $$y'_n(\beta_i)\geq q_{n-K-1}\geq Cq_n.$$
     If $b_{n+1}\beta_i)\neq b_{n+1}(\beta_j)$, then $\beta_i$ and $\beta_j$ are not in the same column in the $n$-towers. Because of the previous results on $x_n$ and $x'_n$, each  of them is at a distance greater than $C\alpha_n$ from the sides of their column, thus 
     $x_n(\beta_i,\beta_j)\geq C\alpha_n$. Otherwise, $x_n(\beta_i,\beta_j)=x_m(\beta_i,\beta_j)$ for the first $m>n$ for which $b_{m+1}(\beta_i)\neq b_{m+1}(\beta_j)$, and we know $m\leq n+K$. Thus for all $n$, 
   $$x_n(\beta_i)\geq C\alpha_{n+K}\geq C\alpha_n.$$
 If $b_{n}\beta_i)\neq b_n(\beta_j)$, then by construction of the towers $y_n(\beta_i,\beta_j))\geq q_{n-1}$. Otherwise, $y_n(\beta_i;\beta_j)=y_m(\beta_i,\beta_j)$ for the last $m<n$ for which $b_{m}(\beta_i)\neq b_m(\beta_j)$, and we know $m\leq n-K$. Thus we get that for all $n$, 
   $$y_n(\beta_i,\beta_j)\geq q_{n-K-1}\geq Cq_n.$$

A cylinder $H$ of length $h$ (for the coding with marked points) is an interval $[y,x[$ for which each  iterate by $R^{-m}$, $1\leq m\leq h-1$, is in a single atom of $Z'$. For a given measure $\mu(H)=x-y$, the minimal value of $h$ is reached when either $x$ and $R_{-h+1}y$, or $y$ and $R_{-h+1}x$, are endpoints of atoms of $Z'$ (otherwise the interval $[y,x[$ could be extended to the left or to the right). We take $n$ such that  $\mu(H)$ is smaller than $\alpha_n$; then as in the proof of Lemma \ref{traj} we see $H$ in the $n$-towers or less than $\alpha_n$ from the right or left of Figure 7. Then the above computations imply that $h$ is at least $Cq_n$ and $\mu(H)$ at least $C\alpha_n$. Hence we get the linear recurrence from Proposition \ref{bos}.\qed\\

The following lemma will be used later.

\begin{lemma}\label{bslr}
If $(X',S)$ is linearly recurrent, when $W$ is a bispecial  word in $L(S)$, of length greater than an initial constant $C_0$, then  if $WU$ is in $L(S)$ with  fixed $\ab U\ab \leq C\ab W\ab$, then $U$ can only be one of two words $U_1$ and $U_2$, where the first letters of $U_1$ and $U_2$  are different, possibly the second letters of $U_1$ and $U_2$  are different, and then the $l$-th letters of $U_1$ and $U_2$ are the same for $l\leq \ab U_1\ab \wedge \ab U_2\ab$. \end{lemma}
{\bf Proof}\\
This is proved by looking in the towers for $S$, using the fact that, by the proof of Theorem \ref{ostlr}, all $y_n(\beta_i)$, $y'_n(\beta_i)$, $y"_n(\beta_i)$ and $y_n(\beta_i,\beta_j)$ are at least $Cq_n$. Then $W$ corresponds to a set of trajectories which coincides on $\ab W\ab$ consecutive symbols, but some (in particular, the leftmost and rightmost ones)  are different on the letter before and the letter after. If all these trajectories are at a distance between $\alpha_{n+1}$ and $\alpha_n$ 
for some $n\geq 2$, then $W$ can be seen in the $n$-towers. 

As $W$ is right special,  it must end just before we see either a $\beta_i$, $i\neq t$,  or $1-\alpha$ between the leftmost and rightmost trajectories in $W$  (as in Lemma \ref{traj} the rightmost one is allowed to hit the considered $\beta_i$ or $1-\alpha$ but not the leftmost). In the first case, these two trajectories disagree on the level containing $\beta_i$; in the second case, the two trajectories disagree left and right of  $1-\alpha$, and on the next letter as they are left and right of $0$; in both cases, then they agree again until we see again $1-\alpha$ or some $\beta_j$ between the leftmost and rightmost trajectories in $W$, thus for a length at least $Cq_n$. As $W$ is left special,  it begins just after we see either a $\beta_i$, $i\neq t$,  or $0$ between the leftmost and rightmost trajectories in $W$, thus by Corollary \ref{ctraj} its length is at most $q_n+q_{n+1}+q_{n+2}\leq C'q_n$, and thus the claimed property is proved. \qed\\

\section{Rigidity for generalized Veech -Sataev}\label{rvs}

\subsection{The natural coding of $T$}
We look now at the {\em natural coding} of $T$, namely its coding by the partition into the $d(r+1)$ intervals used In Definition \ref{dvg} (though they are not necessarily the intervals of continuity of $T$, see Figure 1 above), and we call it $(Y,T)$. We denote by $s_i$ the $i$-th interval in the $s$-th copy of $[0,1[$.
A trajectory $x$ 
 of $T$ under this  natural coding projects on a trajectory $\phi(x)$ of the rotation with marked points $(X',S)$, by applying the map $\phi(s_i)=i$ letter to letter. Because all the $\sigma_i$ are bijective, and their compositions also, as in  Lemma 5 of \cite{fh2} for any word $w$ in $L(T)$, there are exactly $d$ words $v$ such that $\phi (w)=\phi (v)$, and for each of these words either $v=w$ or on the letters $v_i\neq w_i$ for all $i$.\\

As $(X',S)$ is generated by the $r+2$ towers in Figure 7, $(Y,T)$ is generated  by $d(r+2)$ Rokhlin towers. More precisely, by construction of the towers, for all $n$, the trajectories  of the natural coding of $R$
are covered by disjoint occurrences of $M_n$ and $P_n$, the names of the large and small $n$-towers. The trajectories of the coding  with marked points $S$ are covered by the names of the towers in Figure 7: these are denoted by $P_{n,i}$, $1\leq i\leq r_1< r+2$, and $M_{n,j}$, $r_1+1\leq j\leq r+2$, $r_1$ depending on $n$ (we number them from right to left if $n$ is odd, from left to right otherwise). The trajectories of $T$  are covered by  $d(r+2)$ words $P_{n,i,j}$ and $M_{n,i,j}$, $1\leq j\leq d$ which are all the words which project on on $P_{n,i}$ and $M_{n,i}$ by $\phi$.

\begin{pp}\label{vlr} $(Y,T)$ is linearly recurrent if and only if $(X',S)$ is linearly recurrent. In this case, $(Y,T)$ is uniquely ergodic.  \end{pp}
{\bf Proof}\\
Let $[w]$ be a cylinder for $(Y,T)$: for the Lebesgue measure $\mu$ on both sets we have $\mu[w]=\frac1d\mu[\phi w]$, and, for any invariant measure $\nu$ on $(Y,T)$, on $(X',S)$ $\nu$ projects on $\mu$, the unique invariant measure, thus $\mu[\phi w]= \sum_{\phi v=\phi w}\nu [v]\geq \nu[w]$. Hence the result on linear recurrence in both directions comes from Proposition \ref{bos}, while unique ergodicity comes from Proposition \ref{boue}.\qed\\

\subsection{The non-exotic cases}
We use now all the preliminary work to derive results generalizing those in \cite{fh2}, We  do consider these generalizations as non-trivial but do not claim them to be unexpected.

\begin{pp}
If $\alpha$ has unbounded partial quotients,
 $(Y,T)$ is rigid for any invariant measure.
  \end{pp}
{\bf Proof}\\
 In trajectories of $R$, by construction of the towers we have $P_{n+1}=P_n^{a_{n+1}}M_n$, $M_{n+1}=P_n$ for all $n$. Thus $P_{n+2}=(P_n^{a_{n+1}}M_n)^{a_{n+2}}P_n$, $M_{n+2}=P_{n+1}=P_n^{a_{n+1}}M_n$. As $M_n$ is shorter than $P_n$ disjoint occurrences of the word $P_n^{a_{n+1}}$ fill a proportion at least $1-\frac{2}{a_{n+1}+1}$ of the length of both $M_{n+2}$ and $P_{n+2}$. 

In trajectories of $S$, the construction of the towers and the above remark imply that  a proportion at least $1-\frac{2}{a_{n+1}+1}$ of the length of all $M_{n+2,j}$ and $P_{n+2,j}$ is covered by concatenations of the type $P_{n,j_1}...P_{n,j_{a_{n+1}}}$ of length $q_na_{n+1}$. Moreover, all these concatenation contain, at the same place, cycles of the form $P_{n,i_j}^{c_{n,j}}$, where the $c_{n,j}$, $1\leq j\leq r_2\leq r+1$ ($r_2$ depending on $n$) are the successive numbers of $n$-columns containing no $\beta_l$, between two column containing at least one $\beta_l$ or between the sides of the towers and a column containing at least one $\beta_l$ (here column $(0)$ is replaced by its intersection with the large tower). Thus  $\sum_{j=1}^{r_2}c_{n,j}\geq a_{n+1}-r$.

 In trajectories of $T$, we look at the words which project by $\phi$ on cycles $P_{n,i}^c$. 
$n$ and $i$ being fixed, each $P_{n,i,j}$ can be followed by exactly one $P_{n,i,j}$, and thus the $P_{n,i,j}$, $1\leq j\leq d$, are grouped into at most $d$ disjoint strings, each one containing at most $d$ words $P_{n,i,j}$. After the last $P_{n,i,j}$ of  each string, the only $P_{n,i,j'}$ we can see is the first one of the same string. Then, if we move by $T^{d!q_n}$ inside one of the words which project on the cycle
$P_{n,i_j}^{c_{n,j}}$, we go to  the same level in the same tower of name $P_{n,i,j}$, except if we are in the last $d!$ words projecting on this cycle. In each concatenation $P_{n,j_1}...P_{n,j_{a_{n+1}}}$ mentioned above, these ``good" words represent $\sum_{j=0}^{r_2}(c_{n,j}-d!)\vee 0\geq a_{n+1}-r-(r+1)d!$
of the words in $L(T)$ projecting on that concatenation, and thus  a proportion at least $0\vee (1-\frac{2rd!}{a_{n+1}})$ of the length of all $M_{n+2,i,j}$ and $P_{n+2,i,j}$.

 All the levels of the same $n$-tower have the same measure by a given invariant $\mu$, thus if $E$ is a union of levels of the towers of name $P_{n,i,j}$, we have $\mu (E\Delta T^{d!q_n}E)\leq \frac{2rd}{a_{n+1}}$. Now for every set $E$ and $n$ large enough, $E$ can be $\delta_n$-approximated (for the invariant measure $\mu$) by unions of levels of the towers of name $P_{n,i,j}$ or $M_{n,i,j}$; but the towers of name $M_{n,i,j}$ have total measure at most $\frac2{a_{n+1}}$ since they represent a smaller fraction of the length of all $M_{n+2,i,j}$ and $P_{n+2,i,j}$. Thus $\mu (E\Delta T^{d!q_n}E)\leq \frac{2rd+4}{a_{n+1}}+\delta_n$. Hence if the $a_n$ are unbounded $T$ is rigid. \qed\\

The notion of {\em average $\bar d$-separation} is defined in \cite{fh2}, where  comments and explanations on this and related notions can be found.

\begin{definition}\label{asep} For two words of equal length $w=w_1\ldots w_N$ and $w'=w'_1\ldots w'_N$, their Hamming or ${\bar d}$-distance is 
${\bar d}(w,w')=\frac{1}{N}\#\{i; w_i\neq w'_i\}$.\\
A language $L$ on an alphabet $\mathcal A$ is {\em average $\bar d$- separated} for an integer $e\geq 1$ if there exists a language $L'$ on an alphabet $\mathcal A'$, a $K$ to one  (for some $K\geq e$) map $\phi$ from $\mathcal A$ to $\mathcal A'$, extended by concatenation to a map $\phi$ from $L$ to $L'$, such that for any word $w$ in $L$, there are exactly $K$ words $v$ such that $\phi (w)=\phi (v)$, and for each of these words either $v=w$ or $\bar d(w,v)=1$, and a constant $C$, such that 
 if $v_i$ and $v'_i$, $1\leq i\leq e$, are words  in $L$, of equal length $N$, satisfying
\begin{itemize}
 \item 
$\sum_{i=1}^e \bar d(v_i,v'_i)< C$,
\item $\phi (v_i)$ is the same word $u$ for all $i$,
\item $\phi (v'_i)$ is the same word $u'$ for all $i$,
\item $v_i\neq v_j$ for $i\neq j$.
\end{itemize}
Then $\{1,\ldots N\}$ is the disjoint union of three (possibly empty) integer intervals $I_1$, $J_1$, $I_2$ (in increasing order) such that 
\begin{itemize}
\item
$v_{i,J_1}=v'_{i,J_1}$ for all $i$, 
\item $\sum_{i=1}^e\bar d (v_{i,I_1},v'_{i,I_1})\geq 1$ if $I_1$ is nonempty, 
\item $\sum_{i=1}^e\bar d (v_{i,I_2},v'_{i,I_2})\geq 1$ if $I_2$ is nonempty,
\end{itemize}
where $w_{i,H}$ denotes the word made with the $h$-th letters of the word $w_i$ for all $h$ in $H$. \\
This implies in particular that $\# J_1\geq N(1- \sum_{i=1}^e \bar d(v_i,v'_i))$.\\
We call {\em $\bar d$-separation} the average $\bar d$-separation with $K=e=1$, $L=L'$, $\phi$ the identity. \end{definition}

The proof of next proposition will follow step by step the proof of Proposition 44 of \cite{fh2}. The main difference is that in \cite{fh2} $L(T)$ projects by $\phi$ on $L(R)$, while here it projects on the more complicated $L(S)$.  Hence Lemma \ref{bslr} above will replace  Lemma 42 of \cite{fh2}.

 \begin{pp}\label{pdsep}
If the product inequality before Definition \ref{dvg} and the minimality condition of Proposition \ref{pmn} are satisfied,  and the rotation with marked points $(X',S)$ is linearly recurrent, 
$L(T)$ is average $\bar d$-separated with $e=d$.
 \end{pp}
 {\bf Proof}\\
We take $L'=L(S)$, $K=d$. 
Let $v_i$ and $v'_i$   be as   in  Definition \ref{asep}.\\

We compare first $u$ and $u'$;  note that if we see $l$ in some word $\phi (z)$ we see
some $s_l$ at the same place on $z$; 
 thus $\bar d (z,z')\geq \bar d(\phi(z),\phi(z'))$ for all $z$, $z'$; in particular, if $
 \bar d(u,u')=1$, then $\bar d (v_i, v'_i)=1$ for all $i$ and our assertion is proved.
 
 Thus we can assume $\bar d(u,u')<1$.
 We partition $\{1,\ldots N\}$ into
successive integer intervals where $u$ and $u'$ agree or disagree: we get intervals $I_1$, $J_1$,
\ldots, $I_g$, $J_g$, $I_{g+1}$, where $g$ is at least $1$, the intervals are nonempty except 
possibly for $I_1$ or $I_{g+1}$, or both, and for all $j$, $u_{J_j}=u'_{J_j}$, and, except if $I_j$ is empty, $u_{I_j}$ and $u'_{I_j}$ are completely different, i.e.
their distance $\bar d$ is one.

Then for $i\leq g-1$, the word $u_{J_i}=u'_{J_i}$ is right special in the language $L(S)$, and this word is left special if $i\geq 2$.\\

{\em $(H0)$ We suppose first that $u_{J_1}=u'_{J_1}$ is also left special and $u_{J_g}=u'_{J_g}$ is also right special}.\\

Then, by Lemma \ref{bslr}, either $\# J_j$ is smaller than a fixed $m_1$, , or 
 $1 \leq \#I_{j+1}\leq 2$ and $$\#I_{j+1}+
\#J_{j+1}>\ C \#J_j,$$

Similar considerations for $S^{-1}$ imply that for $j>1$ either $\# J_j< m_1$, or $1\leq \#I_{j}\leq 2$ and 
$\#J_{j-1}+\# I_j>C \#J_j$.

We look now at the words $v_i$ and $v'_i$ for some $i$; by the remark above, $v_{i,I_j}$ and $v'_{i,I_j}$ are completely different if $I_j$ is nonempty.
As for $v_{i,J_j}$ and $v'_{i,J_j}$, they have the same image by $\phi$, thus are equal if they begin by the same letter, 
completely different otherwise.\\

Moreover, suppose that $J_j$ has length at least $m_1$, and $v_{i,J_j}=v'_{i,J_j}=Y(i)$, projecting on a right special word $Y$ in $L(S)$ ending with
the letter $j$; then $Y(i)$ ends with the letter $s(i)_j$. Bispecial words in $L(S)$ are described in the proof of Lemma \ref{bslr}; if $Y$ ends just before we see a $\beta_i$, $i\neq t$, after $Y$ in $L(S)$ we see the letters $j_1$ or $j_2$, these denoting the two adjacent intervals around $\beta_i$, then the same $j_3$. Thus after $Y(i)$ in $L(T)$ we see the letters 
$(\sigma s(i))_{j_1}$ or $(\sigma s(i))_{j_2}$ for some $\sigma$, then the letters $(\sigma_{j_4}\sigma s(i))_{j_3}$ or $(\sigma_{j_5}\sigma s(i))_{j_3}$, these two permutations denoting the $\sigma(x)$ on the two mentioned intervals. If $Y$ ends just before we see $1-\alpha$, after $Y$ in $L(S)$ we see the letters $j_6$ or $j_7$, these denoting the two adjacent intervals around $1-\alpha$, then the letters $r$ or $0$, denoting the two adjacent intervals around $0$, then the same $j_{8}$. Thus after $Y(j)$ in $L(T)$, if $\beta_t=1-\alpha$ we see the letters 
$(\sigma s(i))_{j_6}$ and $(\sigma s(i))_{j_7}$ for some $\sigma$, then the letters $(\sigma_{t-1}\sigma s(i))_{0}$ and $(\sigma_{t}\sigma s(i))_{r}$, then the letters $(\sigma_0\sigma_{t-1}\sigma s(i))_{j_{8}}$ and  $(\sigma_r\sigma_{t}\sigma s(i))_{j_{8}}$; if $\beta_j\neq 1-\alpha$ for all $j$, we see the letters $(\sigma s(i))_{j_6}$ and $(\sigma s(i))_{j_7}$ for some $\sigma$, then the letters $(\sigma'\sigma s(i))_{0}$ and $(\sigma'\sigma s(i))_{r}$, for some $\sigma'$, then the letters $(\sigma_0\sigma'\sigma s(i))_{j_{8}}$ and  $(\sigma_r\sigma'\sigma s(i))_{j_{8}}$.
In both cases, this gives us the first letters of the two words $v_{i,J_{j+1}}$ and $v'_{i,J_{j+1}}$.\\

We estimate $c=\sum_{i=1}^d \bar d(v_i,v'_i)$, by looking at 
the indices in some set 
$G_j=
J_j\cup I_{j+1}\cup J_{j+1}$, for any $1\leq j\leq g-1$;
\begin{itemize}
\item if both $\#J_j$ and $\#J_{j+1}$ are
smaller than $m_1$ the contribution of $G_j$ to the sum $c$ is at least $\frac{1}{2m_1+1}$ as
$I_{j+1}$ is nonempty by construction;
\item if $\#J_j\geq m_1$, 
 and 
for at least one $i$ $v_{i,J_j}$ and $v'_{i,J_j}$ are completely different, then the contribution of $G_j$ to $c$
is bigger than $\min(\frac{1}{2}, \frac{C_1}{C_1+1})$ as either $\# J_{j+1}<m_1$ or 
$\# J_j+\# I_{j+1}>C_1\# J_{j+1}$;

\item if 
$\#J_j\geq m_1$ and
for all $i$, $v_{i,J_j}=v'_{i,J_j}=Y(i)$; then, because the $v_i$ are all
different and project by $\phi$ on the same word,  the
 $s(i)$ in the last letter of $Y(i)$
takes $d$ different values when $i$ varies; the hypotheses imply that, in the notations of the previous paragraph
$\sigma_{j_4}\sigma s(i) \neq \sigma_{j_5}\sigma s(i)$ for at least one $i$, and 
$\sigma_0\sigma_{t-1}\sigma s(i))\neq \sigma_r\sigma_{t}\sigma s(i)$, resp. $\sigma_0\sigma'\sigma s(i))\neq \sigma_r\sigma'\sigma s(i)$,  for at least one $i$. 
 This ensures that for this $i$, $v_{i,J_{j+1}}$ and $v'_{i,J_{j+1}}$
 are completely different. As $\#J_{j+1}+\#I_{j+1}>C_1\#J_j$, the contribution of $G_j$ to $c$ 
is bigger than  $C$;

\item if $\#J_{j+1}\geq m_1$, we imitate the last two items by looking in the other direction.

\end{itemize}

Now, if $g$ is even, we can cover $\{1,\ldots Q\}$ by sets $G_j$ and some intermediate $i_l$, and get 
that $c$ is at least a constant $C$. If $g$ is odd and at least $3$,
by deleting either $I_1$ and $J_1$, or $J_g$ and $I_{g+1}$, we cover at least half of $\{1,\ldots Q\}$ 
by sets $G_j$ and some intermediate $i_l$, and $c$ is at least $C$.\\

Thus if $\sum_{i=1}^d \bar d(v_i,v'_i)$ is smaller than a constant $C$, we must have $g=1$; then if 
$\sum_{i=1}^d \bar d(v_i,v'_i)<1$, $v_{i,J_{1}}=v'_{i,J_{1}}$. Thus if $c$ is smaller than  $C$,  we get our conclusion  under the extra hypothesis $(H0)$.\\

For the end of the proof of $\bar d$-separation, without the hypothesis $(H0)$, we refer the reader to the end of the proof of Proposition 44 in \cite{fh2}, as there is nothing different. \qed\\

As is proved in Theorem 3 of \cite{fh2}, for uniquely ergodic systems average $\bar d$-separation implies non-rigidity, but we shall not use that here, as Theorem \ref{nrig} below gives a simpler and more general proof.  We use now the stronger notion of  $\bar d$-separation:

\begin{pp} If the minimality condition is satisfied, $(X',S)$ is linearly recurrent, and,  for all $1\leq u\leq d$,
$\sigma_j(u)\neq \sigma_{j+1}(u)$, $0\leq j \leq r-1$, $j\neq t$,   $\sigma_r\sigma_{t-1} (u)\neq \sigma_0\sigma_t(u)$ (resp. $\sigma_r (u)\neq \sigma_0 (u)$ if $\beta_j\neq 1-\alpha$ for all $j$), then $L(T)$ is $\bar d$-separated and $(Y,T)$ is not of rank one.\end{pp}
{\bf Proof}\\
Then in Proposition \ref{pdsep} we can replace $e=d$ by $e=1$, with the same proof. Then the proof of Theorem 10 of \cite{fh2} applies without modifications. \qed\\

  This last proposition is satisfied in particular for Veech 1969.

\subsection{Isolated points}
In Theorem \ref{ostlr}, the absence of linear recurrence comes from the fact that in the Rokhlin towers some $\beta_i$, $i\neq t$, comes ``too close", horizontally and vertically, to $\alpha$, or $0$ (which is close to $\alpha$), or another $\beta_j$. Thus we define an opposite notion, which we call the isolation (in the towers) of these points.

\begin{definition}
In the Ostrowski expansion of Proposition \ref{dostr}, for integers $n\geq 1$, $M\geq 1$
\begin{itemize} 
\item for $j\neq t$, {\em $\beta_j$ is $(n,M)$-isolated} if there exist  $n-M\leq m_1\leq n$, $n-M\leq m_2\leq n$, $m_2$ odd, $n-M\leq m_3\leq n$, $m_3$ even,  $n-M\leq m'_i\leq n$, $1\leq i\leq r$, $i\neq j$, satisfying
   $b_{m_1}(\beta_j)\neq a_{m_1}-1$, $b_{m_2}(\beta_j)\neq a_{m_2}$, $b_{m_3}(\beta_j)\neq a_{m_3}$, $b_{m'_i}(\beta_j)\neq b_{m'_i}(\beta_i)$, for all
   $1\leq i\leq r$, $i\neq j$,
\item {\em $\alpha$ is $(n,M)$-isolated} if  for all  $1\leq i\leq r$, there exist $n-M\leq m_i\leq n$, $n-M\leq m'_i\leq n$, $m'_i$ odd, $n-M\leq m"_i\leq n$, $m"_i$ even, satisfying
$b_{m_i}(\beta_i)\neq a_{m_i}-1$,
         $b_{m'_i}(\beta_i)\neq a_{m'_i}$, $b_{m"_i}(\beta_i)\neq a_{m"_i}$.
\end{itemize}
\end{definition}

To make statements simpler,  we shall write sometimes that {\em always one of the $\beta_i$ is isolated} to denote there exists  $M$ such that 
 for all $m$, there exists $1\leq j\leq r$, $j\neq t$ such that $\beta_j$ is $(m,M)$-isolated, and, mutatis mutandis, we shall write that {\em always one of the  $\beta_i$ or $\alpha$ is isolated}. \\

Theorem \ref{ostlr} says that $T$ is linearly recurrent whenever always all the $\beta_i$, $i\neq t$, and $\alpha$, are isolated; but weaker assumptions can also be useful. The first one is a sufficient condition for unique ergodicity.  

\begin{pp}\label{cue} If  the minimality condition is satisfied,, $\alpha$ has bounded partial quotients,  and there exists  $\bar M$  and a sequence $p_k\to +\infty$ such that 
  for all $k$, all $\beta_i$, $i\neq t$,  and $\alpha$ are  $(p_k,\bar M)$-isolated, $T$ is uniquely ergodic. \end{pp}
  {\bf Proof}\\ By the proof of Theorem \ref{ostlr}, this condition implies that $\limsup_{n\to +\infty}ne_n(S,\mu)>0$ for the Lebesgue (and unique invariant) measure $\mu$; thus by the proof of Proposition \ref{vlr} $\limsup_{n\to +\infty}ne_n(T,\mu)>0$ for the Lebesgue measure $\mu$, and we conclude by Proposition \ref{boue}. \qed\\
  
Proposition \ref{cue} will be used to build uniquely ergodic examples when needed;   note that its hypothesis  is not equivalent to the unique  ergodicity of $T$: it is not satisfied by Veech 1969 if the sequence $b_n(\beta)$ is made of strings of increasing lengths where either $b_n=a_n-1$ or $b_n=a_n$, $b_{n+1}=0$, while $T$ is uniquely ergodic by \cite{ve69}. \\

Then, as we mentioned above, for our systems non-rigidity will be implied by weaker conditions than the ones ensuring average $\bar d$-separation. 

\begin{theorem}\label{nrig}
Suppose   the minimality condition is satisfied, $\alpha$ has bounded partial quotients,  and there exists  $M$ such that 
\begin{itemize} \item either for all $m$, there exists $1\leq j\leq r$, $j\neq t$, such that $\beta_j$ is $(n,M)$-isolated,
\item or the product inequality is satisfied, and  for all $m$, either $\alpha$ is $(n,M)$-isolated or there exists $1\leq j\leq r$, $j\neq t$, such that $\beta_j$ is $(n,M)$-isolated, \end{itemize}
then $(Y,T)$ is not rigid for any ergodic invariant measure.
\end{theorem}
{\bf Proof}\\
Let $\mu$ be an ergodic invariant measure for $T$. 
Assume that $(Y,T)$ is rigid; then there exists a sequence $Q_k$ tending to infinity such that 
$\mu (D\Delta T^{Q_k}D)$ tends to zero for each of the $d(r+2)$ intervals $D$ defining the natural coding of $T$. 
We fix $\epsilon$ and 
$k$ such that for all these intervals
$$\mu (D\Delta T^{Q_k}D)<\epsilon.$$

Let $A_{D,k}=D\Delta T^{Q_k}D$; 
by the ergodic theorem, for each $D$ and $k$ $\frac{1}{N}\sum_{j=0}^{N-1} 1_{T^jA_{D,k}}(z)$ tends to $\mu (A_{D,k})$, for almost all $z$; we can choose  a set $\Lambda$ of full  $\mu$-measure on which this convergence holds for all $D$ and $k$. Thus for all $z$ in $\Lambda$ and all $k$, there exists $N_0(k)$ such that for all $N$ larger 
than  $N_0(k)$ and 
all $D$, 
$$\frac{1}{N}\sum_{j=0}^{N-1} 1_{T^jA_{D,k}}(z)<\epsilon.$$ By summing these $d(r+2)$ inequalities, we get that 
$$\bar d(z_0\ldots z_{N-1}, z_{Q_k}\ldots z_{Q_k+N-1})<d(r+2)\epsilon$$ for all $N>N_0(k)$. Moreover, for $z$ in some set $\Lambda'$ of full $\mu$-measure, we can
choose $N_0(k)$ 
such that for all $N>N_0(k)$ these inequalities are also satisfied if we replace $z$ by any of the $d$ different points $z'$ such that $\phi (z')=\phi (z)$.\\

We shall now show that this is not possible by estimating $\sum_{i=1}^d\bar d(x^i_0\ldots x^i_{N-1}, y^i_{0}\ldots y^i_{N-1})$ for the $d$ points $x^i$ such that $\phi(x^i)$ is a given point $x$ and the $d$ points $y^i$ such that $\phi(y^i)$ is a given point  $y$. We take $n\geq 1$ such that $\alpha_{n+1}\leq x-y\leq \alpha_n$, and $N$ much larger than $q_n$; we shall look at the trajectories of $x$ and $y$ in the $n$-towers. 

 We partition $\{0,\ldots N-1\}$ into
successive integer intervals where $x$ and $y$ agree or disagree: we get intervals $I_1$, $J_1$,
\ldots, $I_s$, $J_h$, $I_{h+1}$ as in the proof of Proposition \ref{pdsep}; for all $l$, $x_{J_l}=y_{J_l}$, $x_{I_l}$ and $y_{I_l}$ are either empty or completely different, i.e.
their distance $\bar d$ is one.  Except maybe the first one, each $J_l$ begins after we see $\alpha$ or a $\beta_i$, $i\neq t$, between the trajectories of $x$ and $y$, and ends before we see $1-\alpha$ or a $\beta_i$, $i\neq t$, between the trajectories of $x$ and $y$.\\

Suppose that for some $j\neq t$ $\beta_j$ is $(n,M)$-isolated. We group the $I_l$ and $J_l$ into intervals $K_g= I_{l_-(g)}\cup J_{l_-(g)}\cup I_{l_-(g)+1}\cup J_{l_-(g)+1)} ...\cup I_{l_+(g)}\cup J_{l_+(g)}$ where $J_{l_-(g)}$ begins after $\beta_j$, $J_{l_+(g)}$ ends before $\beta_j$, and no other $J_l$ inside $K_g$ has any of these two properties. By  Corollary \ref{ctraj}, for all $g$ $\# K_g\leq 2(q_n+q_{n+1}+q_{n+2})\leq C_1q_n$,
while $\# K_g \geq q_n$ because two times where $\beta_j$ is between the trajectories of $x$ and $y$ are separated by at least $q_n$. 
Also, by the proof  of Theorem \ref{ostlr}, $y_n(\beta_j)$, $y'_n(\beta_j)$, $y"_n(\beta_j)$ and all $y_n(\beta_i,\beta_j)$, $i\neq j$, are at least $C_2q_n$, thus for each $g$ we have $\#J_{l_-(g)}\geq C_2q_n$, $\#J_{l_+(g)}\geq C_2q_n$, 
Now, for each $i$, $x^i_{J_{l_+(g)}}$ and $y^i_{J_{l_+(g)}}$ are either equal or completely different. If for at least one $i$ they are completely different, then $\sum_{i=1}^d\bar d(x^i_{J_{l_+(g)}},y^i_{J_{l_+(g)}})\geq 1$ and $\sum_{i=1}^d\bar d(x^i_{K_g},y^i_{K_g})\geq C_3$. Otherwise, we deduce the first letters of $x^i_{J_{l_-(g+1)}}$ and $y^i_{J_{l_-(g+1)}}$ from the common last letter of 
$x^i_{J_{l_+(g)}}$ and $y^i_{J_{l_+(g)}}$as in the proof of Proposition \ref{pdsep} above, and find  that they must be different for at least one $i$, because the permutations $\sigma_{j_1}$ which is $\sigma(x)$ on the interval left of $\beta_j$ and $\sigma_{j_2}$ on the interval right of $\beta_j$ have different values on at least one point. Then  $\sum_{i=1}^d\bar d(x^i_{J_{l_-(g+1)}},y^i_{J_{l_-(g+1)}})\geq 1$ and $\sum_{i=1}^d\bar d(x^i_{K_{g+1}},y^i_{K_{g+1}})\geq C_3$. Thus we have always $\sum_{i=1}^d\bar d(x^i_{K_g\cup K_{g+1}},y^i_{K_g\cup K_{g+1})}\geq C_4$. We extend $\{0,...N-1\}$ by at most $C_1q_n$ on the left and on the right to a set $K'$ made with an even number of $K_g$; then $\sum_{i=1}^d\bar d(x^i_{K'},y^i_{K'})\geq C_5$ and $\sum_{i=1}^d\bar d(x^i_0\ldots x^i_{N-1}, y^i_{0}\ldots y^i_{N-1})\geq C_5-\frac{2C_1q_n}{N}$.\\

Suppose that $\alpha$ is $(n,M)$-isolated: then we make a similar reasoning. Now our interval $K_g$ are defined by $J_{l_-(g)}$ begins after $\alpha$, $J_{l_+(g)}$ ends before $1-\alpha$, and no other $J_l$ inside $K_g$ has any of these two properties. To get that the first letters of some  $x^i_{J_{l_-(g+1)}}$ and $y^i_{J_{l_-(g+1)}}$ must be different, we use that $\sigma_0\sigma_{t-1}$ and $\sigma_r\sigma_{t}$, resp. $\sigma_0$ and $\sigma_r$ if $\beta_j\neq 1-\alpha$ for all $j$,
have different values on at least one point. By the proof  of Theorem \ref{ostlr}, all $y_n(\beta_ji$, $y'_n(\beta_i)$, $y"_n(\beta_i)$, $1\leq i\leq r$,  are at least $C_2q_n$.
And we get again $\sum_{i=1}^d\bar d(x^i_0\ldots x^i_{N-1}, y^i_{0}\ldots y^i_{N-1})\geq C_5-\frac{2C_1q_n}{N}$.

Under the hypotheses of the theorem, this last relations holds for all $n$ and all $x$ and $y$ with $\alpha_{n+1}\leq x-y\leq \alpha_n$, thus this contradicts rigidity. \qed\\

Theorem \ref{nrig} applies in particular when $(X',S)$ is linearly recurrent, even when $\sigma_0\sigma_t=\sigma_r\sigma_{t-1}$, resp. $\sigma_0=\sigma_r$ if $\beta_j\neq 1-\alpha$ for all $j$ (as soon as there is at least one $\beta_j\neq 1-\alpha$, otherwise we are in the cases of  \cite{fh2}).  As mentioned in the introduction, this gives the first known examples of  non rigid non linearly recurrent interval exchange transformations (note that we could get further examples for any $r\geq 2$ and $d\geq 2$):\\

{\bf Proof of Theorem \ref{nrnlr}}\\ Suppose the conditions of Theorem  \ref{nrig} are satisfied but not those of Theorem \ref{ostlr}. This is possible for example if we build $\beta_1\neq 1-\alpha$ and $\beta_2\neq 1-\alpha$ with prescribed Ostrowski expansions such that, for a fixed  $M$, for all $m$ $\beta_2$ is $(m,M)$-isolated, while there are unbounded strings of consecutive $b_n(\beta_1)=a_n-1$. Then $(X',S)$ is not linearly recurrent and $(Y,T)$ is not rigid, and not linearly recurrent by Proposition \ref{vlr}. Unique ergodicity will be satisfied by Proposition \ref{cue}  if  we ensure   $\beta_1$ is $(p_k,M)$-isolated for a sequence $p_k$. 

 Now, if   we take $r=2$ and $d=2$, with $\beta_1$ and $\beta_2$ as above, and we alternate between the two possible permutations, the identity and the exchange, changing when we cross $\beta_1$ and $\beta_2$, we get the examples claimed in the theorem. \qed\\

When $(X',S)$ is not linearly recurrent, Lemma \ref{bslr} is not satisfied, and we do not know whether $(Y,T)$ is average $\bar d$-separated.

\subsection{In the grey zone: rigidity}
We call grey zone the cases when $\alpha$ has bounded partial quotients, but $(X',S)$ is not linearly recurrent. We could conclude to non-rigidity when   the hypotheses of Theorem \ref{nrig} are  satisfied, but there are still many other  cases. When Theorem \ref{nrig} does not apply, then for all $n$ some $\beta_i$ and $\beta_j$ and/or $\beta_i$, $i\neq t$, and $\alpha$ are too close in the $n$-towers. The simplest case is when all the $\beta_i$ come close to $\alpha$ simultaneously. 

\begin{definition}\label{clu} We say that {\em all the $\beta_i$ cluster on $\alpha$} if
 there exist two sequences $m_k$ and $N_k$, tending to infinity, with $m_k+N_k<m_{k+1}$, such that for all
$1\leq i\leq r$, $i\neq t$, we have 
\begin{itemize} \item either $b_n(\beta_i)=a_n-1$ for all $m_k\leq n\leq m_k+N_k$, 
\item or $b_n(\beta_i)=a_n$ for all even $m_k\leq n\leq m_k+N_k$, 
\item or $b_n(\beta_i)=a_n$ for all odd $m_k\leq n\leq m_k+N_k$.
\end{itemize}\end{definition}

We recall that   $T^n(x,s)=(R^nx, \psi_n(x)s)$ where
$$\psi_n(x)=\sigma(R^{n-1}x)...\sigma(x).$$

\begin{lemma}\label{lclu} Suppose that for a given $n$, for all $i\neq t$, either $x_n(\beta_i, \alpha)<\epsilon \alpha_n$ and $y_n(\beta_i)<\epsilon q_n$, 
or $x_n(\beta_i)<\epsilon \alpha_n$ and $y'_n(\beta_i)<\epsilon q_n$, or $x'_n(\beta_i)<\epsilon \alpha_n$ and $y_n(\beta_i)<\epsilon q_n$; for $0\leq h\leq q_n-1$ we call $\tau_{h,n}$ the permutation $\sigma(x_h)$ when if $n$ is odd $x_h$ is the leftmost (resp. if $n$ is even $x_h$ is the rightmost) point of the level $h$ in the large $n$-tower (the basis being level $0$). Then, on a set $\Xi_n$ of measure at least $1-6\epsilon$, whenever $x$ is in level $h$ of the large or small $n$-tower, then
 $$\psi_{q_n}(x)=\theta_{h,n}=\tau_{h-1,n}...\tau_{0,n}\tau_{q_n-1,n}...\tau_{h,n}.$$
\end{lemma}

\begin{figure}
\begin{center}
\begin{tikzpicture}[scale = 1]

\draw(0,0)--(0,8);\draw(0,0)--(13,0);\draw(13,0)--(13,8);\draw(0,8)--(13,8);
\draw(0,8.15)--(0,11);\draw(0,11)--(3,11);\draw(3,8.15)--(3,11);\draw(0,8.15)--(3,8.15);
\draw[dashed](3,0)--(3,8);
\draw[dashed](6,0)--(6,8);\draw[dashed](9,0)--(9,8);\draw[dashed](12,0)--(12,8);

\draw[dashed](-3,0)--(0,0);
\draw[dashed](-3,0)--(-3,2.85);
\draw[dashed](-3,3)--(-3,11);
\draw[dashed](-3,2.85)--(0,2.85);
\draw[dashed](-3,3)--(0,3);
\draw[dashed](-3,11)--(0,11);
\draw[dashed](0,8)--(0,8.15);
\draw(-3,2.65) node[left]{$0$};
\draw(-3,3.1) node[left]{$\alpha$};

\draw(10,0) node[below]{$\alpha$};

\draw(10,-0.1)--(10,0.1);

\draw(0,11) node[above]{$0$};
\draw(3,11) node[above]{$\alpha_n$};
\draw(0,2.65) node[left]{$\alpha_n$};

\draw[dotted](12.7,0)--(12.7,8);
\draw[dotted](9.7,0)--(9.7,8);
\draw[dotted](10.3,0)--(10.3,8);
\draw[dotted](0.3,0)--(0.3,8);
\draw[dotted](0.3,8.15)--(0.3,11);
\draw[dotted](9.7,0.3)--(13,0.3);
\draw[dotted](-0.3,10.7)--(3,10.7);
\draw[dotted](-0.3,10.7)--(-0.3,11);
\draw[dotted](12.7,7.7)--(13,7.7);
\draw[dotted](-3,2.55)--(-2.7,2.55);
\draw[dotted](-2.7,2.55)--(-2.7,2.85);
\draw[dotted](-3,3.3)--(-2.7,3.3);
\draw[dotted](-2.7,3)--(-2.7,3.3);

\draw(13,8) node[right]{$0$};

\end{tikzpicture}

\caption{All $\beta_i$ cluster on $\alpha$}
\end{center}
\end{figure}
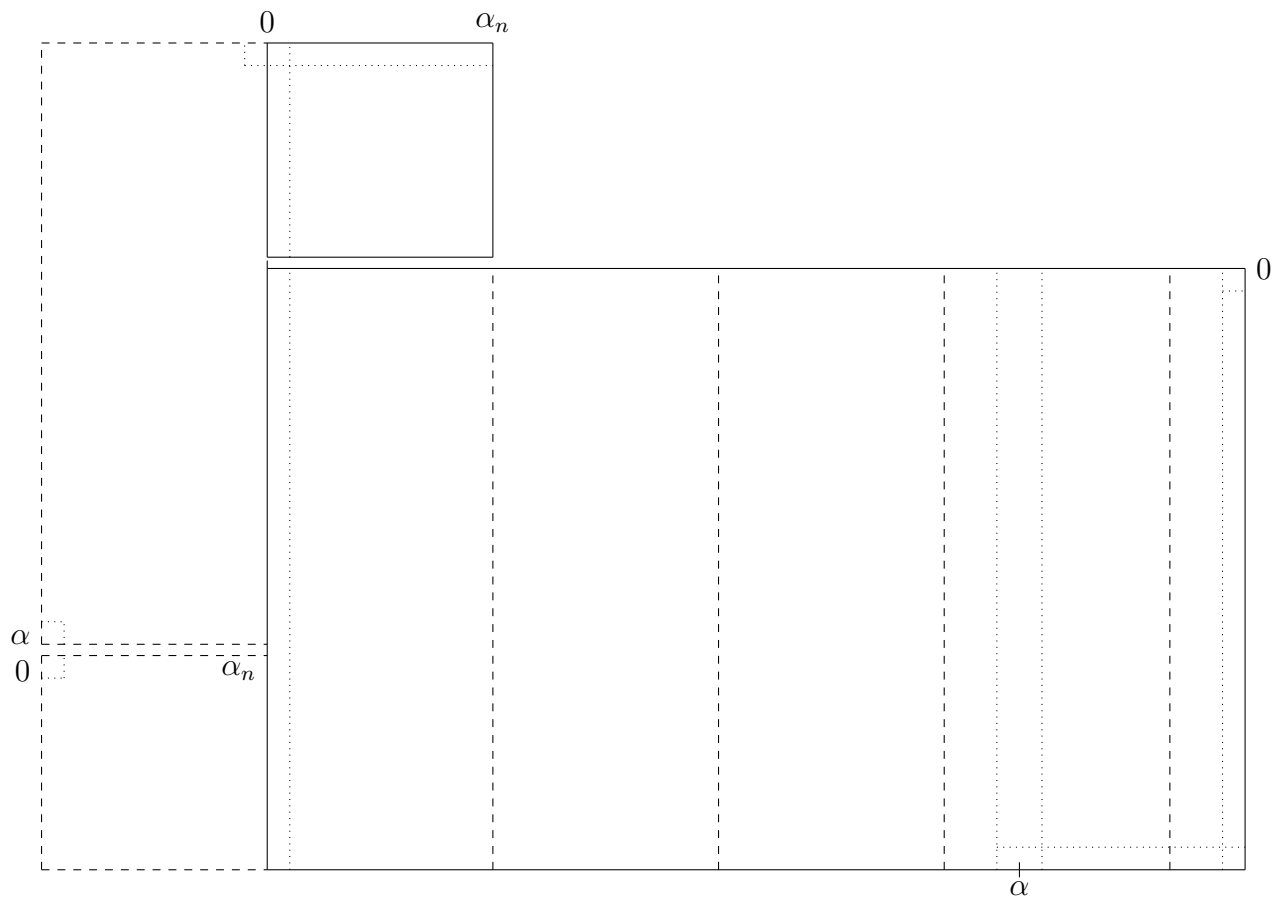

{\bf Proof}\\ We do the proof for $n$ odd.
We delete the set $\Xi_n$, of  small measure as claimed, made with the  $x$  in any of the five following sets:
\begin{itemize} \item the images by $R^m$, $0\leq m\leq q_n-1$, of $[\alpha_n+\alpha -\epsilon\alpha_n, \alpha_n+\alpha[$,
 \item the images by $R^m$, $0\leq m\leq q_n-1$, of $[\alpha -\epsilon\alpha_n, \alpha+\epsilon\alpha_n[$,
  \item the images by $R^m$, $0\leq m\leq \epsilon q_n$, of $[\alpha,\alpha_n+\alpha[$,
  \item the images by $R^m$, $0\leq m\leq q_n+q_{n-1}-1$, of $[\alpha+\alpha_n-\alpha_{n-1}, \alpha+\alpha_n-\alpha_{n-1}+\epsilon\alpha_n[$,
  \item the images by $R^m$, $q_{n-1}-\epsilon q_n\leq m\leq q_{n-1}$, of $[\alpha-\alpha_{n-1}, \alpha - \alpha_{n-1}+\alpha_n[$.
  \end{itemize}
  
 In Figure 7, we show the $n$-towers and what we see less than $\alpha_n$ to their left. The set we delete is between dotted lines, or between dotted line and sides, in the $n$-towers; the $\beta_i$ are confined to the small rectangles near $\alpha$ and $0$ (remember $1-\alpha$ is just below $0$).\\

If $x$ is  in the large $n$-tower but not in column $(0)$, using the above exclusions, we see that whenever the orbit of $x$ is in level $g$, there is no $\beta_i$, $1-\alpha$ or $0$ between this trajectory and $x_g$, and thus the contribution of level $g$ to $\psi_{q_n}(x)$ is $\tau_{g,n}$, and the claimed formula holds.

If $x$ is in the $q_{n-1}$ first levels in column $(0)$, the contribution of level $g$ is $\tau_{g,n}$ until we reach the top of the large $n$-tower; then the orbit of $x$ crosses levels $0$, $1$, ... of the small $n$-tower, staying to the right of $\beta_i$, $1-\alpha$ or $0$. Hence there is no $\beta_i$, $1-\alpha$ or $0$ between this right part  of level $g$ of the small $n$-tower and and the left part of level $g$ of the large $n$-tower,   the contribution of this  level $g$ to $\psi_{q_n}(x)$ is $\tau_{g,n}$, and our formula holds.

If $x$ is in column $(0)$ (either in the large or in the small $n$-tower) above the $q_{n-1}$ first levels but below the upper $\epsilon q_n$ levels (of this column, that is of the small $n$-tower), we continue the reasoning of the previous paragraph. The contributions are the expected ones until we reach the top of the small $n$-tower, whose contribution is $\tau_{q_{n-1}-1,n}$. Then the orbit of $x$ crosses levels $0$, $1$, ... of the large $n$-tower, staying to the right of $\beta_i$, $1-\alpha$ or $0$ (because we have excluded that $x$ is in the leftmost part of width $\epsilon\alpha_n$ of column $(0)$), until we reach the $q_n-1$-th iterate of $x$, which is still at least $\epsilon q_n$ levels below the top. As long as $g\leq q_n-q_{n-1}-1$, there is no $\beta_i$, $1-\alpha$ or $0$ between the right part  of level $g$ of the large $n$-tower and  the left part of level $g+q_{n-1}$ of the large $n$-tower, thus the contributions are the expected ones until the orbit of $x$ reaches level $q_n-q_{n-1}-1$ of the large $n$-tower, whose contribution is $\tau_{q_n-1}$. Then for  $g\geq q_n+q_{n-1}$, there is no $\beta_i$, $1-\alpha$ or $0$ between this right part  of level $g$ of the large $n$-tower and  the left part of level $g+q_{n-1}-q_n$ of the small $n$-tower, there is no $\beta_i$, $1-\alpha$ or $0$ in this level $g+q_n-q_{n-1}$ of the small $n$-tower because the orbit of $x$t has not reached the upper $\epsilon q_n$ levels, there is no $\beta_i$, $1-\alpha$ or $0$ between  this level $g+q_n-q_{n-1}$ of the small $n$-tower and  the left part of the same level of the large $n$-tower. Thus the contributions are still as expected and our result holds. \qed\\

\begin{theorem}\label{tcom} If $\alpha$ has bounded partial quotients, all the $\beta_i$ cluster on $\alpha$, $\sigma_k\sigma_j=\sigma_j\sigma_k$ for all $j,k$, then $(Y,T)$ is rigid for any invariant measure. 
\end{theorem}
{\bf Proof}\\
For any $k$, we choose $n=m_k+[\frac{N_k}2]$. For a given $\epsilon$, by the proof of Theorem \ref{ostlr}, if $k$ is large enough the hypotheses of Lemma \ref{lclu} are satisfied, and its results hold with $\Xi_n$ and $\theta_{h,n}$. 

By definition, for all $h$ and $h'$, $\theta_{h',n}$ is of the form $\theta'\theta_{h,n}\theta'^{-1}$ where $\theta'$ is some composition of the $\sigma(x)$. As all these commute, $\theta_{h,n}$ is a constant $\theta_n$ for all $h$, and $\psi_{q_n}(x)=\theta_{n}$ for all $x$ in $\Xi_n$.

Moreover, if $x$ is in level $h$ in the $n$-towers, $R^{q_n}x$ is in level $h-q_{n-1}$ if $x$ is in column $(0)$ between levels $q_{n-1}$ and $q_n+q_{n-1}-1$, in level   $h+q_n-q_{n-1}$ if $x$ is in the small tower, in level $h$ if $x$ is in any other level. Thus, again as the $\sigma_(x)$ commute, $\psi_{q_n}(R^{q_n}x)=\theta_{n}$ for all $x$ in $R^{-q_n}\Xi_n$, and similarly $\psi_{q_n}(R^{lq_n}x)=\theta_{n}$ for all $x$ in $R^{-lq_n}\Xi_n$, hence $\psi_{lq_n}(x)=\theta_{n}^l$ for all $x$ in $\cap_{l'=0}^{l-1}R^{-l'q_n}\Xi_n$. Let $1\leq \zeta_n\leq d!$ be the order of the permutation $\theta_n$: then $\psi_{\zeta_nq_n}(x)$ is the identity for $x$ in a set of measure at least $1-6d!\epsilon$.

As also $\ab R^{\zeta_nq_n}x -x\ab < \frac{Cd!}{q_n}$, we get that the sequence $\zeta_nq_n$ is a rigidity sequence for $(Y,T)$.\qed\\

Theorem \ref{tcom} is valid in particular when the permutations $\sigma_i$ correspond to the addition of some elements of $\GZ/d\GZ$, as in \cite{ve69}, \cite{ste} or \cite{mer}. The same technique, with more work, applies when the $\sigma_i$ do not commute, but only in some very particular cases.

\begin{pp}\label{3r} Suppose $d=3$, $r=1$, with one marked point $\beta\neq 1-\alpha$, and the two values of $\sigma(x)$ are a transposition and a circular permutation.  For every $\alpha$ with bounded partial quotients, we can find $\beta$ such that $T$ is rigid for any invariant measure.\end{pp}
{\bf Proof}\\
For this, we use again the  quantity $\psi_n(x)=\sigma(R^{n-1}x)...\sigma(x)$.
For any of our systems with $r=1$, we can define by recursion three quantities: 
\begin{itemize} \item 
if $\beta$ is in the large 
$n$-tower, and  $n$ is odd, $\psi_{1,n}$, resp. $\psi_{2,n}$, is the value of  $\psi_{q_n}$ on the basis of the large $n$-tower left, resp. right, of the vertical of $\beta$, $\psi_{3,n}$ is the value of $\psi_{q_{n-1}}$ on the basis of the small $n$-tower;
\item 
if $\beta$ is in the large 
$n$-tower, and  $n$ is even, $\psi_{1,n}$, resp. $\psi_{2,n}$, is the value of  $\psi_{q_n}$ on the basis of the large $n$-tower right, resp. leftt, of the vertical of $\beta$, $\psi_{3,n}$ is the value of $\psi_{q_{n-1}}$ on the basis of the small $n$-tower;
\item 
if $\beta$ is in the small 
$n$-tower, and  $n$ is odd, $\psi_{2,n}$, resp. $\psi_{3,n}$, is the value of  $\psi_{q_{n-1}}$ on the basis of the small $n$-tower left, resp. right, of the vertical of $\beta$, $\psi_{1,n}$ is the value of $\psi_{q_n}$ on the basis of the large $n$-tower;
\item 
if $\beta$ is in the small 
$n$-tower, and  $n$ is even, $\psi_{2,n}$, resp. $\psi_{3,n}$, is the value of  $\psi_{q_{n-1}}$ on the basis of the small $n$-tower right, resp. left, of the vertical of $\beta$, $\psi_{1,n}$ is the value of $\psi_{q_n}$ on the basis of the large $n$-tower;
\end{itemize}

The construction of the towers implies that 
\begin{itemize}
\item if $\beta$ is in the large $n$-tower and $b_{n+1}(\beta)\neq a_{n+1}$, \\
$\psi_{1,n+1}=\psi_{3,n}\psi_{1,n}^{b_{n+1}}\psi_{2,n}^{a_{n+1}-b_{n+1}}$,\\
$\psi_{2,n+1}=\psi_{3,n}\psi_{1,n}^{b_{n+1}+1}\psi_{2,n}^{a_{n+1}-b_{n+1}-1}$,\\
$\psi_{3,n+1}=\psi_{2,n}$;
\item if $\beta$ is in the large $n$-tower and $b_{n+1}(\beta)= a_{n+1}$, \\
$\psi_{1,n+1}=\psi_{3,n}\psi_{1,n}^{a_{n+1}}=\psi_{3,n}\psi_{1,n}^{b_{n+1}}\psi_{2,n}^{a_{n+1}-b_{n+1}}$,\\
$\psi_{2,n+1}=\psi_{2,n}$,\\
$\psi_{3,n+1}=\psi_{1,n}$;
\item if $\beta$ is in the small $n$-tower,\\ 
$\psi_{1,n+1}=\psi_{3,n}\psi_{1,n}^{a_{n+1}}$,\\
$\psi_{2,n+1}=\psi_{2,n}\psi_{1,n}^{a_{n+1}}$,\\
$\psi_{3,n+1}=\psi_{1,n}$.
\end{itemize}

Given $\alpha$, we shall build a $\beta$ clustering on $\alpha$, such that for infinitely many $n$ with $\beta$ close to $\alpha$ in the $n$-tower both $\psi_{1,n}$ and $\psi_{1,n-1}$ are circular permutations, or equivalently have signature $1$.\\

We build $\beta$ by its $b_n(\beta)$, We put $N_0=0$ and  choose an $M_0>N_0+2$; for $N_0+1\leq n\leq M_0-1$, we choose any $0\leq b_n(\beta)\leq a_n-1$, so that  $\beta$ stays in the large $n$-tower, which implies in particular, because of the hypothesis on $\sigma(x)$ and the definition of $T$,  that $\psi_1(n)$ and $\psi_2(n)$ have opposite signatures.
 If $\psi_{1,M_0-1}$ has signature $+1$, we put $M'_0=M_0-1$. Otherwise,  we choose $0\leq b_{M_0}(\beta)\leq a_{M_0}-1$; then if $\psi_{1,M_0}$ has signature $+1$, we put $M'_0=M_0$. If both $\psi_{1,M_0}$ and $\psi_{1,M_0-1}$ have signature $-1$, $\psi_{3,M_0}=\psi_{2,M_0-1}$ and $\psi_{2,M_0}$ have signature $+1$, and the signature of $\psi_{1,M_0+1}$ is $(-1)^{b_{M_0+1}}$; if we choose $b_{M_0+1}$ even this will be $+1$. We choose an even $b_{M_0+1}< a_{M_0+1}$
(this is always possible as we may take $b_{M_0+1}=0$), and put $M'_0=M_0+1$. 

Thus in all cases $\psi_{1,M'_0}$ has signature $+1$ and $\beta$ is in the large $M'_0$-tower. If $\psi_{1,M'_0-1}$ has also signature $+1$, we define $N'_0=M'_0$. Otherwise, $\psi_{1,M'_0-1}$ has  signature $-1$, $\psi_{3,M'_0}=\psi_{2,M'_0-1}$ has signature $+1$, $\psi_{2,M'_0}$ has signature $-1$, and the signature of $\phi_{1,M_0+1}$ is $(-1)^{a_{M_0+1}-b_{M_0+1}}$, and we choose $b_{M_1+1}(\beta)$ so that $\psi_{1,M'_0+1}$ has signature $+1$; if $a_{M_0+1}>1$, we can do it such that $\beta$ is in the large $M_0+1$-tower and put $N'_0=M_0+1$. If $a_{M_0+1}=1$, we choose $b_{M_0+1}=1$ and $\beta$ is in the small $M_0+1$-tower. Using the recursion formulas above, we get that $\psi_{2,M'_0+1}$ has signature $-1$, $\psi_{3,M'_0+1}$ has signature $+1$, $\psi_{1,M'_0+2}$ has signature $+1$, $\psi_{3,M'_0+2}$ has signature $+1$, and $\beta$ is in the large $M'_0+2$-tower. We put $N'_0=M_0+2$.

Then we choose $N_1>N'_0$ and for $N'_0\leq n\leq N_1$ we choose $b_n(\beta)=a_n-1$. The recursion formulas imply that, for all those $n$, $\psi_{1,n}$ has signature $+1$, $\psi_{2,n}$ has signature $-1$, $\psi_{3,n}$ has signature $-1$. Then we choose $b_{N_1+1}\neq a_{N_1+1}-1$ (which may imply that $\beta$ is in the small $N_1+1$-tower), and  start the same process again with $N_0$ replaced by $N_1$. Thus we define sequences $N_k\leq M_k\leq M'_k\leq N'_k<N_{k+1}$, and we choose $N_{k+1}$ so that $N_{k+1}-N'_n$ tends to infinity. \\

We can now adapt the proof of Theorem \ref{tcom}. For any $k$, we choose $n=[\frac{N_{k+1}-N'_k}2]$. For a given $\epsilon$, by the proof of Theorem \ref{ostlr}, if $k$ is large enough the hypotheses of Lemma \ref{lclu} are satisfied, and its results hold with $\Xi_n$ and $\theta_{h,n}$. Moreover, by the proof of Lemma \ref{clu}, $\theta_{0,n}=\psi_{1,n}$. All this is still true if we replace $n$ by $n-1$.

By definition, for all $h$, $\theta_{h',n}$ is of the form $\theta'\theta_{0,n}\theta'^{-1}$ where $\theta'$ is some composition of the $\sigma(x)$. Thus the signature of  $\theta_{h,n}$ is $+1$ for all $h$,  and so is the signature of  $\theta_{h,n-1}$. This implies that all the $\theta_{h,n}$ and  $\theta_{h,n-1}$ are circular permutations, and thus commute. 

We conclude as in Theorem \ref{tcom} that $\psi_{3q_n}(x)$ is some $\theta_{h,n}^3$, and thus the identity,  for all $x$ in a set of measure at least $1-C\epsilon$, and that the sequence $3q_n$ is a rigidity sequence for $(Y,T)$.\qed\\

\subsection{The cases  $d=2$}\label{d2}
In these cases, which constitute the most immediate generalizations of Veech 1969, the two possible permutations  are  the identity $I$ and the exchange $E$. Not only they commute, but, if we have two sequences of such permutations  $\sigma_{i,l}\neq \sigma_{i,r}$ for all $1\leq i\leq K$, then $\sigma_{K,l}...\sigma_{1,l}$ and $\sigma_{K,r}...\sigma_{1,r}$ are equal if $K$ is even, different if $K$ is odd. Thus, as we shall see in the two following propositions, {\em a cluster of an even number of marked points (different from $1-\alpha$) behaves as if there was no marked point at all, and  a cluster of an odd number of such marked points behaves as an isolated marked point}. In theory, with both these properties together with Theorems \ref{nrig} and \ref{tcom},  we could solve completely the question of rigidity for $d=2$ and any number of marked points. However, as the reader may be convinced by studying  Proposition \ref{d2r4} below, a full result would be unduly complicated to state, let alone to prove, so we shall limit ourselves to a complete study of the cases when $1\leq r\leq 3$, and of some examples for $r=4$. These examples in Proposition \ref{d2r4} provide {\em non-rigid examples which do not satisfy the hypotheses of Theorem \ref{nrig}}.\\

\begin{pp} If $\alpha$ has bounded partial quotients and $T$ satisfies the minimality condition, for $d=2$ and at most three marked points different from $1-\alpha$, either Theorem \ref{nrig} applies or $(Y,T)$ is rigid for any invariant measure. 

More precisely: 
\begin{itemize}
\item if $r=1$ and $\beta\neq 1-\alpha$, (the Veech 1969 case), or $r=2$ and $\beta_t=1-\alpha$, $(Y,T)$ is  non-rigid if $(X',S)$ is  linearly recurrent, rigid for any invariant measure otherwise; 
\item if $r=2$ and $\beta_j\neq 1-\alpha$ for all $j$,  or $r=3$ and $\beta_t =1-\alpha$, $(Y,T)$ is non-rigid for any ergodic invariant measure if always one of the $\beta_i$ is  isolated, rigid for any invariant measure otherwise; 
\item if $r=3$ and $\beta_j\neq 1-\alpha$ for all $j$,  or $r=4$ and $\beta_t =1-\alpha$, $(Y,T)$ is non-rigid for any ergodic invariant measure  if always $\alpha$ or one of the $\beta_i$ is  isolated, rigid for any invariant measure otherwise. \end{itemize}
\end{pp}
{\bf Proof}\\
For Veech 1969, Theorem \ref{nrig} does not apply  if and only if $(X',S)$ is not linearly recurrent, and then we can use Theorem \ref{tcom} to get rigidity. This is true also when $r=2$ and $\beta_t=1-\alpha$, as we change permutation, from $I$ to $E$ or from $E$ to $I$,  when we cross $\beta_i$, thus $\sigma_0=\sigma_r$,  $\sigma_t\neq \sigma_{t-1}$, and the product inequality is  satisfied.

When  $r=2$ and $\beta_j\neq 1-\alpha$ for all $j$, $\sigma_0=\sigma_r$;  when  $r=3$ and $\beta_t=1-\alpha$, $\sigma_t\neq \sigma_{t-1}$,  $\sigma_0\neq \sigma_r$, hence in both these cases the product inequality is not satisfied. Thus Theorem \ref{nrig} applies only when always one of the  $\beta_i$ is  isolated, and Theorem \ref{tcom} applies when all $\beta_i$ cluster on $\alpha$. There remains the case where $\alpha$ is always isolated but $\beta_1$ and $\beta_2$ can be very close. In that case, we choose an $n$ such that  $x_n(\beta_1, \beta_2)<\epsilon \alpha_n$
 and $y_n(\beta_1, \beta_2)<\epsilon q_n$. Suppose for example that $\beta_2$ is higher than $\beta_1$ in the $n$-towers; let $\sigma_{i,l}$, resp. $\sigma{i,r}$, be the permutation $\sigma (x)$ on the left (resp. right) of $\beta_i$ on the same level of the $n$-towers, $i=1,2$, let $\sigma_{j_1}, ..., \sigma_{j_h}$ be the values of $\sigma (x)$ on the successive levels
between $\beta_1$ and $\beta_2$. Then $\sigma_{i,l}\neq \sigma_{i,r}$ for $i=1,2$,  and  thus
$\sigma_{2,l}\sigma_{j_h}, ..., \sigma_{j_1}\sigma_{1,l}=\sigma_{2,r}\sigma_{j_h}, ..., \sigma_{j_1}\sigma_{1,r}$ by the remark at the beginning of Section \ref{d2} and commutation. Hence we can make the same reasoning as in Lemma \ref{lclu}: supposing for example that in the $n$-towers $\beta_2$ is higher than $\beta_1$ and to its right, $\beta_i=R^{h_i}\beta'_i$, $i=1,2$, with $\beta'_i$ in the basis of the large $n$-tower; we delete a small set made with the images by $R^m$, $0\leq m\leq q_n-1$, of $[\beta'_1,\beta'_2[$,
  the images by $R^m$, $h_1\leq m\leq h_2$, of the basis of the large $n$-tower, and the upper two levels of the small $n$-tower. Then for the non-deleted $x$  we get the same formula as in Lemma \ref{lclu}, and, as in Theorem \ref{tcom} we conclude that $2q_n$ is a rigidity sequence for $T$.
	
	When $r=3$ and $\beta_j\neq 1-\alpha$ for all $j$, $\sigma_0\neq \sigma_r$; when $r=4$ and $\beta_t=1-\alpha$, we have   $\sigma_t\neq \sigma_{t-1}$ and  $\sigma_0= \sigma_r$, hence in both these cases the product inequality is always satisfied.
Therefore the only case when we cannot apply Theorem \ref{nrig} or Theorem \ref{tcom} is when $\alpha$ and the $\beta_i$ are never isolated, but the $\beta_i$ do not cluster on $\alpha$; thus infinitely often $\alpha$ is close to one of the $\beta_i$, for example $\beta_3$, while $\beta_1$ and $\beta_2$ are very close. For such an $n$, the reasoning of the last case  applies again, and, by deleting all what we have deleted in this case and all we have deleted in Lemma \ref{lclu}, for the non-deleted $x$  we get the same formula as in Lemma \ref{lclu}, and, as in Theorem \ref{tcom} we conclude that $2q_n$ is a rigidity sequence for $T$.\qed\\

In the examples of the next proposition, one of the $\beta_i$ (to make things simpler, we take always the same one, $\beta_1$) will be close to $\alpha$ infinitely often, while the other three will be always far from $\alpha$ and $\beta_1$ but infinitely often close together. This allows non-rigidity though none of our $\beta_i$ or $\alpha$ is always isolated.

\begin{pp}\label{d2r4} Suppose $d=2$, we have four marked points $\beta_1$, $\beta_2$, $\beta_3$, $\beta_4$ different from $1-\alpha$,  the minimality condition is satisfied, $\alpha$ has bounded partial quotients,
 there exist $M_0$ and two sequences $m_k$ and $N_k$, tending to infinity, with $m_k+N_k<m_{k+1}$, such that 
 \begin{itemize}
 \item for all $k$ and all $m_k\leq n\leq m_k+N_k$, $b_n(\beta_1)=a_n-1$, $b_n(\beta_2)=b_n(\beta_3)=b_n(\beta_4)$,
 \item for all $k$ and all $m_k\leq n\leq m_k+N_k+M_0$, there exists $n-M_0\leq m'_1\leq n$ such that $b_{m'_1}(\beta_2)\neq a_{m'_1}-1$,
  \item for all $k$ and all $m_k\leq n\leq m_k+N_k+M_0$, there exists an even $n-M_0\leq m'_2\leq n$ such that $b_{m'_2}(\beta_2)\neq a_{m'_2}$,
   \item for all $k$ and all $m_k\leq n\leq m_k+N_k+M_0$, there exists an odd  $n-M_0\leq m'_3\leq n$ such that $b_{m'_3}(\beta_2)\neq a_{m'_3}$,
 \item for all $k$ and all $m_k+N_k+M_0\leq n\leq m_{k+1}$, $\beta_2$ is $(n,M_0)$-isolated.
 \end{itemize}
Then $T$ is not rigid for any ergodic invariant measure.
\end{pp}
{\bf Proof}\\
By the proof of Theorem \ref{ostlr}, there exists a fixed constant $C_0$, depending only on the size of the partial quotients of $\alpha$, such that, 
\begin{itemize} \item
for any $\beta$, if  there exist $n-2M_0\leq m'_1\leq n$, $n-2M_0\leq m'_2\leq n$, $n-2M_0\leq m'_3\leq n$ such that $m'_2$ is even, $m'_3$ is odd, $b_{m'_1}(\beta)\neq a_{m'_1}-1$,
$b_{m'_2}(\beta)\neq a_{m'_2}$,
$b_{m'_3}(\beta)\neq a_{m'_3}$,
then both $y_n(\beta)$ and $y'_n(\beta)$ are
 at least $C_0q_n$;
  \item
for any $\beta\neq \beta'$, if  there exists $n-2M_0\leq m'_4\leq n$ such that $b_{m'_4}(\beta)\neq b_{m'_4}(\beta')$,
then  $y_n(\beta,\beta')$ is 
 at least $C_0q_n$. \end{itemize}
 
 Our hypotheses ensure that or our system, the first  result holds  for every $n$ with $\beta=\beta_2$, and also (because of the values of  $b_n(\beta_1)$, $b_n(\beta_2)$, $b_n(\beta_3)$, $b_n(\beta_4)$ for $m_k\leq n\leq m_k+N_k$) that both results hold for $\beta=\beta_2$, $\beta=\beta_3$, $\beta=\beta_4$, $\beta'=\beta_1$ for $m_k+M_0\leq n\leq m_k+N_k+M_0$ (that is why we have chosen $2M_0$ to define $C_0$). 
 
 Using the other part of the proof of Theorem \ref{ostlr}, we choose $M_1>M_0$, depending only on the size of the partial quotients of $\alpha$, such that,
\begin{itemize}
\item for any $\beta\neq \beta'$, if $b_{m'}(\beta)=b_{m'}(\beta')$ for all $n\leq m'\leq n+M_1$, $x_n(\beta, \beta')\leq \frac{\alpha_{n+1}}{4}$ (remember that $\alpha_{n+1}\geq C\alpha_n$),
\item for any $\beta\neq \beta'$, if $b_{m'}(\beta)=b_{m'}(\beta')$ for all $n-M_1\leq m'\leq n$, $y_n(\beta, \beta')\leq \frac{C_0q_n}{2}$,
\item for any $\beta$, if $b_{m'}(\beta)=a_{m'}-1$ for all $n-M_1\leq m'\leq n$, $y_n\leq \frac{C_0q_n}{2}$.
\end{itemize} 

Now we make the beginning of the proof of Theorem \ref{nrig} above: to contradict rigidity, we have  to estimate $\sum_{i=1}^d\bar d(x^i_0\ldots x^i_{N-1}, y^i_{0}\ldots y^i_{N-1})$ for the $d$ points $x^i$ such that $\phi(x^i)$ is a given point $x$ and the $d$ points $y^i$ such that $\phi(y^i)$ is a given point  $y$. We take $n\geq 1$ such that $\alpha_{n+1}\leq \rho=x-y\leq \alpha_n$, and $N$ much larger than $q_n$; we shall look at the trajectories of $x$ and $y$ in the $n$-towers. \\

Suppose $m_k+M_1\leq n\leq m_k+N_k-M_1$. For this $n$, we place $\beta_2$, $\beta_3$, $\beta_4$ in the $n$-towers. We call $\beta$ the one which is lowest, $\beta"$ the highest, $\beta'$ the middle one. As in the proof of Theorem \ref{nrig}, we cut $\{0,...N-1\}$ into intervals
 $I_l$ and $J_l$ and group them  into intervals $K_g= I_{l_-(g)}\cup J_{l_-(g)}\cup I_{l_-(g)+1}\cup J_{l_-(g)+1)} ...\cup I_{l_+(g)}\cup J_{l_+(g)}$ where $J_{l_-(g)}$ begins after $\beta$, $J_{l_+(g)}$ ends before $\beta$, and no other $J_l$ inside $K_g$ has any of these two properties.  We have again that  for all $g$ $\# K_g\leq C_1q_n$. and $\# K_g \geq q_n$.
 
 The beginning of $J_{l_-(g)}$ and the end of $J_{l_+(g)}$ correspond to a $j$ such that $\beta$ is between $T^jx$ and $T^jy$, which is equivalent to $T^jy\in [\beta-\rho, \beta[$; by the ergodic theorem, for $N$ large, there are about $\rho N\geq \alpha_{n+1}N$ such indices $j$. We call ``bad" those $j$ for which $T^jy$ is in $[\beta-\frac{\alpha_{n+1}}{4}, \beta[$ or $T^jy$ is in $[\beta-\rho, \beta-\rho+\frac{\alpha_{n+1}}{4} [$, which correspond at most to about $N\frac{\alpha_{n+1}}{2}$ indices. By deleting all $K_g$ for which 
 $J_{l_+(g)}$  ends before a bad $j$, we keep at least half of the intervals $K_g$. 
 Again, we look at the transition between $K_g$ and $K_{g+1}$ for the non-deleted $K_g$. The beginning  of $J_{l_+(g)}$ is $\alpha$ or a $\beta_i$; the possible one making $J_{l_+(g)}$ shortest is either $\alpha$ or $\beta_1$, which is at least $C_0q_n$ far (vertically) from $\beta$; thus $\#J_{l_+(g)}$ is at least $C_0q_n$.
For each $i$, $x^i_{J_{l_+(g)}}$ and $y^i_{J_{l_+(g)}}$ are either equal or completely different. If for at least one $i$ they are completely different, this gives a contribution of $1$ to the global $\bar d$-sum on the length of $J_{l_+(g)}$.

Now, by our hypothesis,  both $\beta'$ and $\beta"$ are $\frac{\alpha_{n+1}}{4}$ close (horizontally) to $\beta$ and $\frac{C_0q_n}{2}$  close (vertically) to $\beta$. Thus the fact that our $K_g$ has not been deleted  guarantees that after seeing $\beta$ between $T^jy$ and $T^jx$, we shall see  $\beta'$ between $T^{j'}y$ and $T^{j'}x$, $\beta"$ between $T^{j"}y$ and $T^{j"}x$,, with $j<j'<j"<j+\frac{C_0q_n}{2}$; and we do not see either $1-\alpha$ or $\beta_1$ before  as we are far enough from the top of the towers. Thus $J_{l_-(g+1)+2}$ begins with $\beta"$, and ends before a point which, in the case that makes it shortest, is either $\beta_1$ or $1-\alpha$ and is at least $\frac{C_0q_n}2$
far (vertically) from $\beta"$.

If $x^i_{J_{l_+(g)}}$ and $y^i_{J_{l_+(g)}}$ are equal for all $i$, we shall deduce from their common last letter the first letters of $x^i_{J_{l_-(g+1)+2}}$ and $y^i_{J_{l_-(g+1)+2}}$  as in the proof of Proposition \ref{pdsep} above. For that we use again the remark at the beginning of Section \ref{d2}: let $\sigma_{1,l}$, $\sigma_{2,l}$, $\sigma_{3,l}$, resp.  $\sigma_{1,r}$, $\sigma_{2,r}$, $\sigma_{3,r}$ be the permutations $\sigma (x)$ on the left (resp. right) of $\beta$, $\beta'$, $\beta"$ on the same level of the $n$-towers. The two permutations involved in computing the letter we want are, by commutation,  $\sigma\sigma_{3,l}\sigma_{2,l}\sigma_{1,l}$  and  $\sigma\sigma_{3,r}\sigma_{2,r}\sigma_{1,r}$  for a fixed $\sigma$, and these are different. This gives a contribution of $1$ to the global $\bar d$-sum on the length of $J_{l_-(g+1)}+2$.

Thus, for each non-deleted $K_g$, there is a contribution of $1$ to the global $\bar d$-sum on a length at least $\frac{C_0q_n}{2}\geq \frac{C_0}{2C_1}\#K_g$. The non-deleted $K_g$ make a proportion at least $\frac{1}{2C_1}$ of $\{0,...N-1\}$, thus the global $\bar d$-sum cannot be close to $0$.\\

Suppose now  $m_k+N_k+M_0\leq n\leq m_{k+1}+M_1$. Then $\beta_2$ is $(n,M_0+M_1)$ isolated and, after fixing $x$ and $y$  we conclude as in the proof of Theorem \ref{nrig} that the global $\bar d$-sum cannot be close to $0$.\\

Suppose now $m_k+N_k-M_1\leq n\leq m_k+N_k+M_0$. For these $n$, our hypotheses ensure that there exist $n \leq m'_4\leq n+M_1+M_0$ such that $b_{m'_4}(\beta_2)\neq b_{m'_4}(\beta_3)$,
$n \leq m'_5\leq n+M_1+M_0$ such that $b_{m'_5}(\beta_2)\neq b_{m'_5}(\beta_4)$, By the proof of Theorem \ref{ostlr}, this implies that both $x_n(\beta_2,\beta_3)$ and $x_n(\beta_2,\beta_4)$ are at least $C_2\alpha_n$.

We fix an $n$ and place  $\beta_2$, $\beta_3$, $\beta_4$ in the $n$-towers. Again, we fix $x=y+\rho$, define the $I_l$ and $J_l$. 
\begin{itemize} \item If $\beta_2$ is the leftmost of the points $\beta_2$, $\beta_3$, $\beta_4$. By the ergodic theorem, for $N$ large, there are about $\rho N\leq \alpha_{n}N$  indices $j$ such that $\beta_2\leq T^jy\leq \beta_2+\rho$, and at least about $C_2\alpha_{n}N$  indices $j$ such that $\beta_2\leq T^jy\leq \beta_2+C_2\alpha_n$; 
\item if $\beta_2$ is the rightmost of the points $\beta_2$, $\beta_3$, $\beta_4$. By the ergodic theorem, for $N$ large, there are about $\rho N\leq \alpha_{n}N$  indices $j$ such that $\beta_2-\rho\leq T^jx\leq \beta_2$, and at least about $C_2\alpha_{n}N$  indices $j$ such that $\beta_2-C_2\alpha_n\leq T^jx\leq \beta_2$; 
\item If $\beta_2$ is the middele one of the points $\beta_2$, $\beta_3$, $\beta_4$, suppose for example $\beta_3$ is the leftmost one. By the ergodic theorem, for $N$ large, there are about $\rho N\leq \alpha_{n}N$  indices $j$ such that $\beta_3\leq T^jy\leq \beta_3+\rho$, and at least about $C_2\alpha_{n}N$  indices $j$ such that $\beta_3\leq T^jy\leq \beta_3+C_2\alpha_n$.
\end{itemize}

We group the $I_l$ and $J_l$ in intervals $K_g$, using $\beta=\beta_2$ in the first two cases, $\beta=\beta_3$ in the last case. Take the first case for example:   for a proportion at least $C_2$ of the $K_g$,  $J_{l_+(g)}$  ends at a  $j$ such that $T^jx$ is to the right of $\beta_2$, and between $\beta_2$ and the verticals of $\beta_3$ and $\beta_4$. Hence for these $j$ we cannot see $\beta_3$ or $\beta_4$ between the trajectories of $x$ and $y$ before $j$ and after the basis of the towers, or after $j$ and before the top of the towers; thus for these $K_g$ the permutations giving the first letter of $\#J_{l_-(g+1)}$ are the same as when $\beta_2$ is isolated. The vertical distances from $\beta_2$ to $\alpha$, $1-\alpha$ and $\beta_1$ being bounded from below as in the previous case, both $\#J_{l_+(g)}$  and $\#J_{l_-(g+1)}$ are  at least $C_0q_n$; thus  for this proportion $C_2$ of the $K_g$ there is a contribution of $1$ to the global $\bar d$-sum on a length at least $C_0q_n\geq \frac{C_0}{C_1}\#K_g$. The other cases are similar, and we conclude that the global $\bar d$-sum cannot be close to $0$.
 \qed\\

Note that in the particular case of $d=2$, two different permutations are different on all points, so we could make the above reasonings on each $\bar d(x^i_0\ldots x^i_{N-1}, y^i_{0}\ldots y^i_{N-1})$, but that would not simplify significantly the computations. 

We can make examples satisfying the hypotheses of Proposition \ref{d2r4} for every value of $\alpha$. For example, if all $a_n$ are equal to $1$, for $m_k\leq n\leq m_k+N_k$, $b_n(\beta_1)$ will always be $0$ while $b_n(\beta_2)=b_n(\beta_3=b_n(\beta_4)$ can be successively $1,0,0,1,0,0, 1,0,0...$


\begin{thebibliography}{50}


\bibitem{afh} P. ARNOUX, S. FERENCZI, P. HUBERT: Trajectories of rotations, Acta Arith. 87 (1999), no. 3, 209--217.



\bibitem{bo} M. BOSHERNITZAN: Rank two interval exchange transformations, Ergodic Theory Dynam. Systems 8 (1988), no. 3, 379--394.

\bibitem{bo2} M. BOSHERNITZAN: A condition for unique ergodicity of minimal symbolic flows, Ergodic Theory Dynam. Systems 12 (1992), no. 3, 425--428. 



\bibitem{fr1} S. FERENCZI: Systems of finite rank, Colloq. Math. 73 (1997), 35--65. 




\bibitem{fh2} S. FERENCZI, P. HUBERT: Rigidity of interval exchanges, J. Mod. Dyn. 14 (2019), 153--177.



\bibitem{fogg} S. FERENCZI, T. MONTEIL:  Infinite words with uniform frequencies, and invariant measures, Combinatorics, automata and number theory, 373--409, Encyclopedia Math. Appl., 135 (2010), Cambridge Univ. Press, Cambridge.


\bibitem{gp} M. GUENAIS, F. PARREAU: Valeurs propres de transformations li\'ees aux rotations
irrationnelles et aux fonctions en escalier (eigenvalues of transformations arising from irrational
rotations and step functions, (French), preprint, arXiv: 0605250.

\bibitem{ma} H. MASUR: Interval exchange transformations and measured foliations, Annals of Mathematics,
115 (1982), 169--200.

\bibitem{ma2} H. MASUR: Hausdorff dimension of divergent Teichm\" uller geodesics, Trans. Amer. Math. Soc. 324 (1991), no. 1, 235--254.

\bibitem{mat} H. MASUR, S. TABACHNIKOV:
Rational billiards and flat structures, Handbook of dynamical systems, Vol. 1A, 1015--1089, North-Holland, Amsterdam, 2002. 

\bibitem{mer} K. D. MERRILL: Cohomology of step functions under irrational rotations, Isr.
J. of Math. 52 (1985), 320--340.

\bibitem{robe} D. ROBERTSON: Mild mixing of certain interval exchange transformations,  Ergodic Theory Dynam. Systems 39 (2019), no. 1, 248--256.


\bibitem{sat}  E. A. SATAEV: The number of invariant measures for flows on orientable surfaces, (Russian)  Izv. Akad. Nauk SSSR Ser. Mat.  39  (1975), no. 4, 860--878, translated in Mathematics of the USSR-Izvestiya,  9 (1975), 813--830.

\bibitem{ste} M. STEWART: Irregularities of uniform distribution, Acta Math. Acad. Sc.
Hung. 37 (1981), 1--39.

\bibitem{ve69} W. A. VEECH: Strict ergodicity in zero dimensional dynamical systems and the Kronecker-Weyl theorem mod 2, Trans. Amer. Math. Soc. 140, (1969), 1--33.

 \bibitem{vs} W. A. VEECH: A criterion for a process to be prime, Monatsh. Math. 94 (1982), no. 4, 335--341.
 
 \bibitem{vb} W. A. VEECH: Boshernitzan's criterion for unique ergodicity of an interval exchange transformation, Ergodic Theory Dynam. Systems 7 (1987), no. 1, 149--153.



 

 
\end{thebibliography}
\end{document}